# CENTRAL LIMIT THEOREM FOR SIGNAL-TO-INTERFERENCE RATIO OF REDUCED RANK LINEAR RECEIVER

By G. M. Pan[1] and W. Zhou[2]

*EURANDOM and National University of Singapore*

Let $\mathbf{s}_k = \frac{1}{\sqrt{N}}(v_{1k},\ldots,v_{Nk})^T$, with $\{v_{ik}, i, k = 1, \ldots\}$ independent and identically distributed complex random variables. Write $\mathbf{S}_k = (\mathbf{s}_1,\ldots,\mathbf{s}_{k-1},\mathbf{s}_{k+1},\ldots,\mathbf{s}_K)$, $\mathbf{P}_k = \mathrm{diag}(p_1,\ldots,p_{k-1},p_{k+1},\ldots,p_K)$, $\mathbf{R}_k = (\mathbf{S}_k\mathbf{P}_k\mathbf{S}_k^* + \sigma^2\mathbf{I})$ and $\mathbf{A}_{km} = [\mathbf{s}_k, \mathbf{R}_k\mathbf{s}_k, \ldots, \mathbf{R}_k^{m-1}\mathbf{s}_k]$. Define $\beta_{km} = p_k\mathbf{s}_k^*\mathbf{A}_{km}(\mathbf{A}_{km}^* \times \mathbf{R}_k\mathbf{A}_{km})^{-1}\mathbf{A}_{km}^*\mathbf{s}_k$, referred to as the signal-to-interference ratio (SIR) of user $k$ under the multistage Wiener (MSW) receiver in a wireless communication system. It is proved that the output SIR under the MSW and the mutual information statistic under the matched filter (MF) are both asymptotic Gaussian when $N/K \to c > 0$. Moreover, we provide a central limit theorem for linear spectral statistics of eigenvalues and eigenvectors of sample covariance matrices, which is a supplement of Theorem 2 in Bai, Miao and Pan [*Ann. Probab.* **35** (2007) 1532–1572]. And we also improve Theorem 1.1 in Bai and Silverstein [*Ann. Probab.* **32** (2004) 553–605].

## 1. Introduction.

1.1. *The signal-to-interference ratio (SIR) in engineering.* Consider a synchronous direct-sequence code-division multiple-access (CDMA) system. Suppose that there are $K$ users and that the dimension of the signature sequence $\mathbf{s}_k$ assigned to user $k$ is $N$. Let $x_k$ denote the symbol transmitted by user $k$, $p_k$ the power of user $k$ and $\mathbf{n} \in \mathbb{C}^N$ noise vector with mean zero and covariance matrix $\sigma^2\mathbf{I}$. Suppose that $x_k's$ are independent random variables (r.v.'s) with $Ex_k = 0$ and $Ex_k^2 = 1$ and that $x_k's$ are independent of $\mathbf{n}$. The

Received January 2007; revised June 2007.
[1]Supported in part by NSFC Grant 10471135 and 10571001 and by NUS Grant R-155-050-055-133/101.
[2]Supported in part by NUS Grant R-155-000-076-112.

discrete time model for the received vector $\mathbf{r}$ is

$$\mathbf{r} = \sum_{k=1}^{K} \sqrt{p_k} x_k \mathbf{s}_k + \mathbf{n}. \tag{1.1}$$

The goal in wireless communication is to estimate the transmitted $x_k$ for each user in an appropriate receiver. For simplicity, in the sequel we are only interested in linear receivers. A linear receiver, represented by a vector $\mathbf{c}_k$, estimates $x_k$ in a form $\mathbf{c}_k^* \mathbf{r}$ (the notation $^*$ denotes the complex conjugate transpose of a vector or matrix). The well known linear mean-square error (MMSE) minimizes

$$E|x_k - \mathbf{c}_k^* \mathbf{r}|^2. \tag{1.2}$$

To evaluate the linear receivers, a popular performance measure is the output signal-to-interference ratio (SIR),

$$\frac{p_k (\mathbf{c}_k^* \mathbf{s}_k)^2}{\sigma^2 \mathbf{c}_k^* \mathbf{c}_k + \sum_{j \neq k}^{K} p_j (\mathbf{c}_k^* s_j)^2} \tag{1.3}$$

(see Verdú [19] or Tse and Hanly [16]). Ideally, a good receiver should have a higher SIR.

Without loss of generality we focus only user 1. For MMSE receiver, from (1.2) one can solve $\mathbf{c}_1 = \mathbf{R}_1^{-1} \mathbf{s}_1$ and then substitute $\mathbf{c}_1$ into (1.3) to obtain the SIR expression for user 1 as

$$\hat{\beta}_1 = p_1 \mathbf{s}_1^* \mathbf{R}_1^{-1} \mathbf{s}_1, \tag{1.4}$$

where $\mathbf{R}_1 = (\mathbf{S}_1 \mathbf{P}_1 \mathbf{S}_1^* + \sigma^2 \mathbf{I})$, $\mathbf{S}_1 = (\mathbf{s}_2, \ldots, \mathbf{s}_K)$ and $\mathbf{P}_1 = \text{diag}(p_2, \ldots, p_K)$. It turns out that the choice of $\mathbf{c}_1$ also maximizes user 1's SIR. But since MMSE involves a matrix inverse this may be very costly when the spreading factor is high. Based on this reason, some simple and near MMSE performance receivers like reduced-rank linear receiver have been considered.

The basic idea behind a reduced rank is to project the received vector onto a lower dimensional subspace. For the multistage Wiener (MSW), the lower dimensional subspace has been described as a set of recursions by Goldstein, Reed and Scharf [7] and Honig and Xiao [10]. However, we would like to make use of another property of MSW given in Theorem 2 in Honig and Xiao [10] for our purpose, that is, MSW receiver estimates $x_1$ through MMSE after producing $m$-dimensional project vector $\mathbf{A}_{1m}^* \mathbf{r}$ instead of $\mathbf{r}$, where $m < n$ and

$$\mathbf{A}_{1m} = [\mathbf{s}_1, \mathbf{R}_1 \mathbf{s}_1, \ldots, \mathbf{R}_1^{m-1} \mathbf{s}_1]. \tag{1.5}$$

Similar to (1.4), one can get $c_{1m} = (\mathbf{A}_{1m}^* \mathbf{R}_1 \mathbf{A}_{1m})^{-1} \mathbf{A}_{1m}^* \mathbf{s}_1$ and the output SIR

$$\boxed{\beta_{1m} = p_1 \mathbf{s}_1^* \mathbf{A}_{1m} (\mathbf{A}_{1m}^* \mathbf{R}_1 \mathbf{A}_{1m})^{-1} \mathbf{A}_{1m}^* \mathbf{s}_1,} \tag{1.6}$$

which is the focus of this paper.

The MSW, as a kind of reduced-rank receiver, was first introduced by Goldstein, Reed and Scharf [7]. The receiver is widely employed in practice because the number of stages $m$ needed to achieve a target SIR, unlike other reduced-rank receivers, does not scale with the system size, that is, dimensionality $N$ of the system, as remarked by Honig and Xiao [10]. In their subsequent newsletter [11], the authors specially addressed this point. In addition, Honig and Xiao [10] showed that the SIR of MSW converges to a deterministic limit in a large system. However, as we know, in a finite system, the SIR will fluctuate around the limit. Moreover, such fluctuation will lead to some important performance measures, such as error probability and outrage probability. Regarding this promising receiver, we will characterize such fluctuation by providing central limit theorems in this paper.

From now on the signature sequences are modeled as random vectors, that is,

$$\mathbf{s}_k = \frac{1}{\sqrt{N}}(v_{1k}, \ldots, v_{Nk})^T,$$

$k = 1, \ldots, K$, where $\{v_{ik}, i, k = 1, \ldots\}$ are independent and identically distributed (i.i.d.) r.v.'s. Then the SIRs (1.6) may be further analyzed using the random matrices theory when $K$ and $N$ go to infinity with their ratio being a positive constant, which is well known as the large system analysis in the wireless communication field.

Tse and Hanly [16] and Verdú and Shamai [20] derived, respectively, the large system SIR and spectral efficiency under MMSE, Matched filter (MF) and decorrelator receiver. Tse and Zeitouni [17] proved that the distribution of SIR under MMSE is asymptotically Gaussian. Later, Bai and Silverstein [4] reported the asymptotic SIR under MMSE for a general model. For more progress in this area, one may see the review paper of Tulino and Verdú [18] and, in addition, refer to the review paper of Bai [2] concerning random matrices theory. Here we would also like to say a few words about our earlier work (Pan, Guo and Zhou [12]). In that paper, the random variables are assumed to be real and we could apply central limit theorems which have appeared in the literature. For example, we made use of main results from Götze and Tikhomirov ([8], page 426: considering real random variables with the sixth moment) and Bai and Silverstein [3] (requiring $Ev_{11}^4 = 3$ or $E|v_{11}|^4 = 2$). In the present work we develop a central limit theorem for the statistic of eigenvalues and eigenvectors under the finite fourth moment (see Theorem 1.3), which further gives a central limit theorem for a random quadratic form (see Remark 1.5). And we give a central limit theorem (see Theorem 1.4) for eigenvalues by dropping the assumption $Ev_{11}^4 = 3$ or $E|v_{11}|^4 = 2$ in Bai and Silverstein [3]. For central limit theorems in other matrix models, we refer to [1].

Our main contribution to engineering is to prove that the distribution of the SIR under MSW, after scaling, is asymptotic Gaussian and that the sum of the SIRs for all users under MF ($m = 1$), after subtracting a proper value, has a Gaussian limit, which further gives the asymptotic distribution of the sum mutual information under MF.

We introduce some notation before stating our results. Set $\mathbf{R} = (\mathbf{C} + \sigma^2 \mathbf{I})$, $\mathbf{C} = \mathbf{SPS}^*$, $\mathbf{S} = (\mathbf{s}_1, \ldots, \mathbf{s}_K)$ and $\mathbf{P} = \mathrm{diag}(p_1, \ldots, p_K)$. Suppose that $F^{c,H}(x)$ and $H(x)$, respectively, denote the weak limit of the empirical spectral distribution function $F^{c_N \mathbf{SPS}^*}$ and $H_N$ (i.e., $F^{\mathbf{P}}$), where $c_N = N/K$. In particular, $F^{c,H}(x)$ becomes $F^c(x)$ when $\mathbf{P}$ is the identity matrix, whose probability density was given in Jonsson [6]. Let $W^0(t)$ denote a Brownian bridge and $X$ is independent of $W^0(t)$, which is $N(0, Ev_{11}^4 - 1)$. Furthermore, let

$$W_x^c = W^0(F^c(x)),$$

$$\zeta_i = \sum_{u=0}^{i} \binom{i}{u} (\sigma^2)^{i-u} \left( h_u X + \sqrt{2} c^{-u} \int_{(1-\sqrt{c})^2}^{(1+\sqrt{c})^2} x^u \, dW_x^c \right),$$

$i = 1, \ldots, 2m - 1$, and $\zeta_0 = X$, with $h_u = \int x^u \, dF^c(cx)$. Define $a_m = \int (x + \sigma^2)^m \, dF^{c_N, H_N}(cx)$ and

$$\mathbf{b} = (1, a_1, \ldots, a_{m-1})^T, \qquad \mathbf{B} = \begin{pmatrix} a_1 & a_2 & \cdots & a_m \\ a_2 & a_3 & \cdots & a_{m+1} \\ \cdots & \cdots & \cdots & \cdots \\ a_m & a_{m+1} & \cdots & a_{2m-1} \end{pmatrix},$$

where $F^{c_N, H_N}(x) = F^{c,H}(x)|_{c = c_N, H = H_N}$.

In what follows, with a slight abuse of notation, we still use $a_m$ as a limit, such as (1.8) below, even when $F^{c_N, H_N}(x)$ is replaced by $F^{c,H}(x)$ in the expression of $a_m$.

THEOREM 1.1. *Suppose that:*

(a) $\{v_{ij}, i, j = 1, \ldots,\}$ *are i.i.d. complex r.v.'s with* $Ev_{11} = 0$, $Ev_{11}^2 = 0$, $E|v_{11}|^2 = 1$ *and* $E|v_{11}|^4 < \infty$.
(b) $c_N \to c > 0$ *as* $N \to \infty$.
(c) $p_1 = \cdots = p_K = 1$. *Then, for any finite integer* $m$,

(1.7) $$\sqrt{N}(\beta_{1m} - \mathbf{b}^* \mathbf{B}^{-1} \mathbf{b}) \xrightarrow{D} y,$$

*where*

(1.8) $$y = 2\boldsymbol{\zeta}^* \mathbf{B}^{-1} \mathbf{b} - \mathbf{b}^* \mathbf{B}^{-1} \mathbf{D} \mathbf{B}^{-1} \mathbf{b},$$

*with* $\boldsymbol{\zeta}^* = (\zeta_0, \ldots, \zeta_{m-1})$ *and* $\mathbf{D} = (d_{ij}) = (\zeta_{i+j-1})$.

REMARK 1.1. It can be verified that

$$\text{Cov}\left(\int_{(1-\sqrt{c})^2}^{(1+\sqrt{c})^2} x^i \, dW_x^c, \int_{(1-\sqrt{c})^2}^{(1+\sqrt{c})^2} x^j \, dW_x^c\right)$$

(1.9)
$$= \int_{(1-\sqrt{c})^2}^{(1+\sqrt{c})^2} x^{i+j} \, dF^c(x)$$
$$- \int_{(1-\sqrt{c})^2}^{(1+\sqrt{c})^2} x^i \, dF^c(x) \int_{(1-\sqrt{c})^2}^{(1+\sqrt{c})^2} x^j \, dF^c(x).$$

Moreover, $X$ is independent of $\int_{(1-\sqrt{c})^2}^{(1+\sqrt{c})^2} x^i \, dW_x^c$ and so the variance of $y$ can be computed, although it is complicated.

The asymptotic distribution of the sum mutual information has been derived for MMSE by Pan, Guo and Zhou [12]. Thus, it is interesting to derive the corresponding asymptotic distribution of the MSW. But, unfortunately, it is rather complicated for the MSW case. At this stage, we can only derive the asymptotic distribution for the sum mutual information for the case $m = 1$, which is well known as the MF (see Verdú [19]).

Obviously, when $m = 1$, the output SIR for the MSW, $\beta_{km}$ (the expressions for $\beta_{km}$ can be derived similarly to $\beta_{1m}$), becomes

(1.10)
$$\beta_k = \frac{p_k (\mathbf{s}_k^* \mathbf{s}_k)^2}{\mathbf{s}_k^* \mathbf{R}_k \mathbf{s}_k},$$

with $\mathbf{R}_k = \mathbf{C}_k + \sigma^2 \mathbf{I}$ and $\mathbf{C}_k = \mathbf{S}_k \mathbf{P}_k \mathbf{S}_k^*$, where $\mathbf{S}_k$ and $\mathbf{P}_k$ are respectively obtained from $\mathbf{S}$ and $\mathbf{P}$ by deleting the $k$th column (here we denote $\beta_{k1}$ by $\beta_k$).

THEOREM 1.2. *Suppose that:*

(a) $\{v_{ij}, i, j = 1, \ldots, \}$ *are i.i.d. complex r.v.'s. with* $Ev_{11} = 0$, $Ev_{11}^2 = 0$, $E|v_{11}^2| = 1$ *and* $E|v_{11}|^4 < \infty$.
(b) *The empirical distribution function of power matrix* $\mathbf{P}$ *converges weakly to some distribution function* $H(t)$ *with all the powers bounded by some constant.*
(c) $c_N \to c > 0$ *as* $N \to \infty$. *Then*

(1.11)
$$\sum_{k=1}^{K} \left(\beta_k - \frac{p_k}{\sigma^2 + c}\right) \xrightarrow{D} N(\mu, \tau^2)$$

*with* $p_1 = \cdots = p_K = 1$,

$$\mu = \frac{2E|v_{11}|^4 - 3}{c(\sigma^2 + 1/c)^2} + \frac{1}{c^2(\sigma^2 + 1/c)^3},$$

*and* $\tau$ *defined in (5.34).*

We would like to point out that the result has been given only for the equal power case ($p_1 = \cdots = p_K = 1$) in Theorem 1.2, although the assumptions are concerning different powers. As will be seen, the main difficulty of the different powers case is that matrices $(\mathbf{SPS}^*)^2$ and $\mathbf{SP}^2\mathbf{S}^*$ have different eigenvalues. But, it is worth pointing out that one may establish a central limit theorem for

$$\sum_{j=1}^{N}(f(\lambda_j) + g(\mu_j))$$

following a similar line of Bai and Silverstein [3], where $f, g$ are analytical functions and $\lambda_j, \mu_j$ denote the eigenvalues of $\mathbf{P}^{1/2}\mathbf{S}^*\mathbf{SP}^{1/2}$ and $\mathbf{PS}^*\mathbf{SP}$, respectively. We do not intend to pursue this direction since the process is lengthy.

Concerning the sum mutual information under the MF, we have the following:

COROLLARY 1.1. *Under the conditions of Theorem 1.2,*

$$(1.12) \qquad \sum_{k=1}^{K}\left(\log(1+\beta_k) - \log\left(1 + \frac{1}{\sigma^2+c}\right)\right) \xrightarrow{D} N(\mu_1, \tau_1^2)$$

*with*

$$\mu_1 = \frac{\mu}{1+(c^{-1}+\sigma^2)^{-1}}$$
$$- \frac{2(E|v_{11}|^4 - 2)(c^{-1}+\sigma^2)^2 + 2c^{-1}(1+c^{-1}) + \sigma^4 + 2\sigma^2 c^{-1}}{c(c^{-1}+\sigma^2)^4(1+(c^{-1}+\sigma^2)^{-1})^2}$$

*and*

$$\tau_1^2 = \frac{\tau^2}{(1+(c^{-1}+\sigma^2)^{-1})^2}.$$

1.2. *Random matrices.* Random matrices have been used in wireless communication since Grant and Alexander's 1996 conference presentation [9] and it has proved to be a very powerful technique. To prove the preceding theorems, we develop a central limit theorem for the eigenvalues and eigenvectors of the sample covariance matrices, which is a supplement of Theorem 2 in Bai, Miao and Pan [5]. And we also improve Theorem 1.1 in Bai and Silverstein [3]. Obviously, these central limit theorems are interesting themselves.

Let $c_N \mathbf{T}_N^{1/2}\mathbf{SS}^*\mathbf{T}_N^{1/2} = \mathbf{A}_N$ with $\mathbf{T}_N^{1/2}$ being the square root of a nonnegative definite matrix $\mathbf{T}_N$ and $\mathbf{U}_N \Lambda_N \mathbf{U}_N^*$ be the spectral decomposition of $\mathbf{A}_N$, where $\Lambda_N = \text{diag}(\lambda_1, \lambda_2, \ldots, \lambda_N)$, $\mathbf{U}_N = (u_{ij})$ is a unitary matrix consisting of the orthonormal eigenvectors of $\mathbf{A}_N$. Suppose that $\mathbf{x}_N =$ $(x_{N1},\ldots,x_{NN})^T \in \mathbb{C}^N, \|\mathbf{x}_N\| = 1$, is nonrandom and $\mathbf{y} = (y_1, y_2, \ldots, y_N)^T = U_N^* \mathbf{x}_N$. Let $F^{\mathbf{A}_N}$ denote the empirical spectral distribution (ESD) of the matrix $\mathbf{A}_N$ and $F_1^{\mathbf{A}_N}(x)$ another ESD of $\mathbf{A}_N$, that is,

$$(1.13) \qquad F_1^{\mathbf{A}_N}(x) = \sum_{i=1}^{N} |y_i|^2 I(\lambda_i \leq x).$$

Let

$$G_N(x) = \sqrt{N}(F_1^{\mathbf{A}_N}(x) - F^{c_N, H_N}(x)),$$

and $\underline{m}(z) = \underline{m}_{F^{c,H}}(z)$ denote the Stieltjes transform of the limiting empirical distribution function of $c_N \mathbf{S}^* \mathbf{T}_N \mathbf{S}$. Now it is time to state the following theorem.

THEOREM 1.3. *Assume:*

(1) $v_{ij}, i, j = 1, 2, \ldots,$ *are i.i.d. with* $Ev_{11} = 0, E|v_{11}|^2 = 1$ *and* $E|v_{11}|^4 < \infty$, *and* $\lim_{N \to \infty} c_N = c \in (0, \infty)$.
(2) $\mathbf{x}_N \in \{\mathbf{x} \in \mathbb{C}^N, \|\mathbf{x}\| = 1\}$.
(3) $\mathbf{T}_N$ *is nonrandom Hermitian nonnegative definite such that its spectral norm is bounded in* $N$, $H_N = F^{\mathbf{T}_N} \xrightarrow{\mathcal{D}} H$, *a proper distribution function and* $\mathbf{x}_N^*(T - zI)^{-1}\mathbf{x}_N \to m_{F^H}(z)$, *where* $m_{F^H}(z)$ *denotes the Stieltjes transform of* $H(t)$.
(4) $g_1, \ldots, g_k$ *are defined and analytic on an open region* $\mathcal{D}$ *of the complex plane, which contains the real interval*

$$(1.14) \qquad \left[\liminf_N \lambda_{\min}^{\mathbf{T}_N} I_{(0,1)}(c)(1-\sqrt{c})^2, \limsup_N \lambda_{\max}^{\mathbf{T}_N}(1+\sqrt{c})^2\right],$$

*where* $\lambda_{\min}^{\mathbf{T}_N}$ *and* $\lambda_{\max}^{\mathbf{T}_N}$ *denote, respectively, the minimum and maximum eigenvalues of* $\mathbf{T}_N$.

(5)
$$\sup_z \sqrt{N}\left|x_N^*(\underline{m}_{F^{c_N,H_N}}(z)\mathbf{T}_N + I)^{-1}x_N - \int \frac{1}{\underline{m}_{F^{c_N,H_N}}(z)t+1}dH_N(t)\right| \to 0,$$

*as* $n \to \infty$.
(6)
$$\max_i |\mathbf{e}_i^* \mathbf{T}_N^{1/2}(z\underline{m}(z)\mathbf{T}_N + zI)^{-1} x_N| \to 0,$$

*where* $\mathbf{e}_i$ *is the* $N \times 1$ *column vector with the* $i$*th element being* 1 *and the rest being* 0. *Then the following conclusions hold:*

(a) If $v_{11}$ and $\mathbf{T}_N$ are real, the random vector $(\int g_1(x)\,dG_N(x), \ldots, \int g_k(x)\,dG_N(x))$ converges weakly to a Gaussian vector $(X_{g_1}, \ldots, X_{g_k})$, with mean zero and covariance function

$$
\begin{aligned}
&\mathrm{Cov}(X_{g_1}, X_{g_2}) \\
(1.15) \quad &= -\frac{1}{2\pi^2} \int_{\mathcal{C}_1} \int_{\mathcal{C}_2} g_1(z_1) g_2(z_2) \\
&\qquad \times \frac{(z_2 \underline{m}(z_2) - z_1 \underline{m}(z_1))^2}{c^2 z_1 z_2 (z_2 - z_1)(\underline{m}(z_2) - \underline{m}(z_1))}\, dz_1\, dz_2.
\end{aligned}
$$

The contours $\mathcal{C}_1$ and $\mathcal{C}_2$ in the above equality are disjoint, both contained in the analytic region for the functions $(g_1, \ldots, g_k)$ and both enclosing the support of $F^{c_n, H_n}$ for all large $n$.

(b) If $v_{11}$ is complex, with $Ev_{11}^2 = 0$, then the conclusion (a) still holds, but the covariance function reduces to half of the quantity given in (1.15).

REMARK 1.2. It is under the assumption $Ev_{11}^4 = 3$ in the real case or $E|v_{11}|^4 = 2$ in the complex case that Bai, Miao and Pan [5] obtained the above central limit theorem. But, when $Ev_{11}^4 \neq 3$ in the real case, there exist sequences $\{\mathbf{x}_n\}$ such that

$$\left( \int x\,dG_N(x), \int x^2\,dG_N(x) \right)$$

fails to converge in distribution, as pointed out in Silverstein [13]. Therefore, when $Ev_{11}^4 \neq 3$ in the real case or $E|v_{11}|^4 \neq 2$ in the complex case, to guarantee the central limit theorem, we here impose an additional condition (6), which is implied by

$$(1.16) \qquad \max_i |x_{Ni}| \to 0,$$

when $\mathbf{T}_N$ becomes a diagonal matrix. Thus, the variance is dependent on the fourth moment of $v_{11}$.

REMARK 1.3. Let $g_1(x) = x, g_2(x) = x^2, \ldots, g_k(x) = x^k$. Then

$$\sqrt{N}\left( \left( \mathbf{x}_N^* \mathbf{A}_N \mathbf{x}_N - \int x\,dF^{c_n, H_N}(x) \right), \ldots, \left( \mathbf{x}_N^* \mathbf{A}_N^k \mathbf{x}_N - \int x^k\,dF^{c_n, H_N}(x) \right) \right)$$

converges weakly to a Gaussian vector, which is used when proving Theorem 1.1.

To derive Theorem 1.2, we would like to present a central limit theorem for the eigenvalues, which is a little improvement of Theorem 1.1 in Bai and Silverstein [3]. Define

$$L_N(x) = N(F^{\mathbf{A}_N}(x) - F^{c_N, H_N}(x)).$$

THEOREM 1.4. *In addition to the assumptions* (1), (3) *and* (4) *in Theorem 1.3 [remove the assumption concerning* $\mathbf{x}_N^*(\mathbf{T}_N - zI)^{-1}\mathbf{x}_N$ *in* (3)*], suppose that*

$$\frac{1}{N}\sum_{i=1}^{N}\mathbf{e}_i^*\mathbf{T}_N^{1/2}(\underline{m}(z_1)\mathbf{T}_N+\mathbf{I})^{-1}\mathbf{T}_N^{1/2}\mathbf{e}_i\mathbf{e}_i^*\mathbf{T}_N^{1/2}(\underline{m}(z_2)\mathbf{T}_N+\mathbf{I})^{-1}\mathbf{T}_N^{1/2}\mathbf{e}_i$$

(1.17)
$$\to h_1(z_1, z_2)$$

*and*

(1.18)
$$\frac{1}{N}\sum_{i=1}^{N}\mathbf{e}_i^*\mathbf{T}_N^{1/2}(\underline{m}(z)\mathbf{T}_N+\mathbf{I})^{-1}$$
$$\times\mathbf{T}_N^{1/2}\mathbf{e}_i\mathbf{e}_i^*\mathbf{T}_N^{1/2}(\underline{m}(z)\mathbf{T}_N+\mathbf{I})^{-2}\mathbf{T}_N^{1/2}\mathbf{e}_i \to h_2(z).$$

*Then the following conclusions hold:*

(a) *If* $v_{11}$ *and* $\mathbf{T}_N$ *are real, then* $(\int g_1(x)\,dL_N(x),\ldots,\int g_k(x)\,dL_N(x))$ *converges weakly to a Gaussian vector* $(X_{g_1},\ldots,X_{g_k})$, *with mean*

(1.19)
$$EX_g = -\frac{1}{2\pi i}\int g(z)\frac{c\int \underline{m}^3(z)t^2\,dH(t)/(1+t\underline{m}(z))^3}{(1-c\int \underline{m}^2(z)t^2\,dH(t)/(1+t\underline{m}(z))^2)^2}\,dz$$
$$-\frac{Ev_{11}^4 - 3}{2\pi i}\int g(z)\frac{c\underline{m}^3(z)h_2(z)}{1-c\int \underline{m}^2(z)t^2\,dH(t)/(1+t\underline{m}(z))^2}\,dz$$

*and covariance function*

(1.20)
$$\mathrm{Cov}(X_{g_1}, X_{g_2})$$
$$= -\frac{1}{2\pi^2}\int\int \frac{g_1(z_1)g_2(z_2)}{(\underline{m}(z_1)-\underline{m}(z_2))^2}\frac{d}{dz_1}\underline{m}(z_1)\frac{d}{dz_2}\underline{m}(z_2)\,dz_1\,dz_2$$
$$- \frac{c(Ev_{11}^4 - 3)}{4\pi^2}\int\int g_1(z_1)g_2(z_2)\frac{d^2}{dz_1\,dz_2}$$
$$\times [\underline{m}(z_1)\underline{m}(z_2)h_1(z_1,z_2)]\,dz_1\,dz_2.$$

(b) *If* $v_{11}$ *is complex with* $Ev_{11}^2 = 0$, *then* (a) *holds as well, but the mean is now*

(1.21) $$EX_g = -\frac{E|v_{11}|^4 - 2}{2\pi i}\int g(z)\frac{c\underline{m}^3(z)h_2(z)}{1-c\int \underline{m}^2(z)t^2\,dH(t)/(1+t\underline{m}(z))^2}\,dz$$

*and covariance function*

$$\mathrm{Cov}(X_{g_1}, X_{g_2})$$
$$= -\frac{1}{4\pi^2}\int\int \frac{g_1(z_1)g_2(z_2)}{(\underline{m}(z_1)-\underline{m}(z_2))^2}\frac{d}{dz_1}\underline{m}(z_1)\frac{d}{dz_2}\underline{m}(z_2)\,dz_1\,dz_2$$

(1.22)
$$-\frac{c(E|v_{11}|^4 - 2)}{4\pi^2} \int\int g_1(z_1)g_2(z_2)\frac{d^2}{dz_1\,dz_2} \times [\underline{m}(z_1)\underline{m}(z_2)h_1(z_1,z_2)]\,dz_1\,dz_2.$$

REMARK 1.4. When $\mathbf{T}_N$ is a diagonal matrix,
$$h_2(z) = \int \frac{t^2\,dH(t)}{(\underline{m}(z)t+1)^3},$$
$$h_1(z_1,z_2) = \int \frac{t^2\,dH(t)}{(\underline{m}(z_1)t+1)(\underline{m}(z_2)t+1)}.$$

This indicates that the assumptions $Ev_{11}^4 = 3$ or $E|v_{11}|^4 = 2$ in Bai and Silverstein [3] can be removed when $\mathbf{T}_N$ is a diagonal matrix. When $\mathbf{T}_N = \mathbf{I}$ and $g(x) = x^r$,

$$\frac{1}{2\pi i}\int g(z)\frac{c\underline{m}^3(z)h_2(z)}{1 - c\int \underline{m}^2(z)t^2\,dH(t)/(1+t\underline{m}(z))^2}\,dz$$

(1.23)
$$= c^{1+r}\sum_{j=0}^{r}\binom{r}{j}\left(\frac{1-c}{c}\right)^j\binom{2r-j}{r-1}$$
$$- c^{1+r}\sum_{j=0}^{r}\binom{r}{j}\left(\frac{1-c}{c}\right)^j\binom{2r+1-j}{r-1},$$

and when $g_1(x) = x^{r_1}$ and $g_2(x) = x^{r_2}$,

$$-\frac{c}{4\pi^2}\int\int g_1(z_1)g_2(z_2)\frac{d^2}{dz_1\,dz_2}[\underline{m}(z_1)\underline{m}(z_2)h_1(z_1,z_2)]\,dz_1\,dz_2$$

(1.24)
$$= c^{r_1+r_2+1}\sum_{j_1=0}^{r_1}\sum_{j_2=0}^{r_2}\binom{r_1}{j_1}\binom{r_2}{j_2}\left(\frac{1-c}{c}\right)^{j_1+j_2}$$
$$\times \binom{2r_1-j_1}{r_1-1}\binom{2r_2-j_2}{r_2-1}.$$

REMARK 1.5. In applying Theorem 1.4 to Theorem 1.2, we take $g_1(x) = x + x^2$, that is, one needs to transform (1.11) into

$$\sum_{j=1}^{n}(\lambda_j + \lambda_j^2) + u_n,$$

where the term $u_n$ will be proved to converge to some constant in probability. Indeed, when using Theorem 1.3 or Theorem 1.4, $g_1(x)$ is usually taken to be a polynomial function.

The rest of this paper is organized as follows. The proofs of Theorem 1.3 and Theorem 1.1 are given in Sections 2 and 3, respectively. Section 4 includes the argument of Theorem 1.4. Section 5 establishes Theorem 1.2, while the truncation of the underlying r.v.'s is postponed until the Appendix. Section 6 establishes Corollary 1.1. Throughout this paper, to save notation, $M$ may denote different constants on different occasions.

**2. Proof of Theorem 1.3.** Let $\mathbf{A}(z) = \mathbf{A}_N - z\mathbf{I}$, $\mathbf{A}_j(z) = \mathbf{A}(z) - \mathbf{s}_j\mathbf{s}_j^*$. With a slight abuse of notation, here and in the argument of Theorem 1.4, we use $\mathbf{s}_j$ to denote the $j$th column of $c_N^{1/2}\mathbf{T}_N^{1/2}\mathbf{S}$, as in Bai, Miao and Pan [5], but one should note that this $\mathbf{s}_j$ is different from one of other parts. To complete the proof of Theorem 1.3, according to the argument of Theorem 2 in Bai, Miao and Pan [5] [especially (4.1), (4.5) and (4.7)], it is sufficient to prove that

$$(2.1) \quad \frac{1}{K}\sum_{j=1}^{K}\sum_{i=1}^{N} E_j(\mathbf{H}_{nj}(z_1))_{ii} E_j(\mathbf{H}_{nj}(z_2))_{ii} \xrightarrow{i.p.} 0,$$

where $E_j = E(\cdot|\mathcal{F}_j)$, $\mathcal{F}_j = \sigma(\mathbf{s}_1,\ldots,\mathbf{s}_j)$ and

$$\mathbf{H}_{nj}(z) = \mathbf{T}_N^{1/2}\mathbf{A}_j^{-1}(z)\mathbf{x}_n\mathbf{x}_n^*\mathbf{A}_j^{-1}(z)\mathbf{T}_N^{1/2}.$$

Define

$$\mathbf{A}_{jk}(z) = \mathbf{A}(z) - \mathbf{s}_j\mathbf{s}_j^* - \mathbf{s}_k\mathbf{s}_k^*, \qquad \varepsilon_k(z) = \beta_{jk}(z)\mathbf{A}_{jk}^{-1}(z)\mathbf{s}_k\mathbf{s}_k^*\mathbf{A}_{jk}^{-1}(z),$$

$$E\hat{\mathbf{H}}_{nj}(z) = \mathbf{T}_N^{1/2}E\mathbf{A}_j^{-1}(z)\mathbf{x}_n\mathbf{x}_n^*E\mathbf{A}_j^{-1}(z)\mathbf{T}_N^{1/2}, \qquad \beta_{jk}(z) = \frac{1}{1+\mathbf{s}_k^*\mathbf{A}_{jk}(z)\mathbf{s}_k}.$$

It is observed that

$$\mathbf{e}_i^*\mathbf{T}_N^{1/2}(\mathbf{A}_j^{-1}(z_1) - E\mathbf{A}_j^{-1}(z_1))\mathbf{x}_n\mathbf{x}_n^*\mathbf{A}_j^{-1}(z_1)\mathbf{T}_N^{1/2}\mathbf{e}_i$$

$$= \mathbf{e}_i^*\mathbf{T}_N^{1/2}(\mathbf{A}_j^{-1}(z_1) - E\mathbf{A}_j^{-1}(z_1))\mathbf{x}_n\mathbf{x}_n^*(\mathbf{A}_j^{-1}(z_1) - E\mathbf{A}_j^{-1}(z_1))\mathbf{T}_N^{1/2}\mathbf{e}_i$$

$$\quad + \mathbf{e}_i^*\mathbf{T}_N^{1/2}(\mathbf{A}_j^{-1}(z_1) - E\mathbf{A}_j^{-1}(z_1))\mathbf{x}_n\mathbf{x}_n^*E\mathbf{A}_j^{-1}(z_1)\mathbf{T}_N^{1/2}\mathbf{e}_i$$

$$= \sum_{k_1,k_2=1}^{K}\mathbf{e}_i^*\mathbf{T}_N^{1/2}(E_{k_1}\mathbf{A}_j^{-1}(z_1) - E_{k_1-1}\mathbf{A}_j^{-1}(z_1))\mathbf{x}_n$$

$$\qquad \times \mathbf{x}_n^*(E_{k_2}\mathbf{A}_j^{-1}(z_1) - E_{k_2-1}\mathbf{A}_j^{-1}(z_1))\mathbf{T}_N^{1/2}\mathbf{e}_i$$

$$\quad + \sum_{k=1}^{K}\mathbf{e}_i^*\mathbf{T}_N^{1/2}(E_k\mathbf{A}_j^{-1}(z_1) - E_{k-1}\mathbf{A}_j^{-1}(z_1))\mathbf{x}_n\mathbf{x}_n^*E\mathbf{A}_j^{-1}(z_1)\mathbf{T}_N^{1/2}\mathbf{e}_i$$

$$= \sum_{k_1 \neq j, k_2 \neq j}^{K} \mathbf{e}_i^* \mathbf{T}_N^{1/2}((E_{k_1} - E_{k_1-1})\varepsilon_{k_1}(z_1))\mathbf{x}_n$$
$$\times \mathbf{x}_n^*((E_{k_2} - E_{k_2-1})\varepsilon_{k_2}(z_1))\mathbf{T}_N^{1/2}\mathbf{e}_i$$
$$- \sum_{k \neq j}^{K} \mathbf{e}_i^* \mathbf{T}_N^{1/2}((E_k - E_{k-1})\varepsilon_k(z_1))\mathbf{x}_n\mathbf{x}_n^* E\mathbf{A}_j^{-1}(z_1)\mathbf{T}_N^{1/2}\mathbf{e}_i.$$

This, together with the Burkholder inequality and (4.4) in Bai, Miao and Pan [5], gives

$$\left[E\left|\sum_{i=1}^{N} E_j(\mathbf{H}_{nj}(z_1) - \mathbf{T}_N^{1/2}(E\mathbf{A}_j^{-1}(z_1))\mathbf{x}_n\mathbf{x}_n^*\mathbf{A}_j^{-1}(z_1)\mathbf{T}_N^{1/2})_{ii}(E_j\mathbf{H}_{nj}(z_2))_{ii}\right|\right]^2$$

$$\leq \sum_{i=1}^{N} E|(\mathbf{H}_{nj}(z_1) - \mathbf{T}_N^{1/2}(E\mathbf{A}_j^{-1}(z_1))\mathbf{x}_n\mathbf{x}_n^*\mathbf{A}_j^{-1}(z_1)\mathbf{T}_N^{1/2})_{ii}|^2$$

$$\times \sum_{i=1}^{N} E|(\mathbf{H}_{nj}(z_2))_{ii}|^2$$

$$\leq M \sum_{i=1}^{N} \left[E\left|\sum_{k_1 \neq j} \mathbf{e}_i^* \mathbf{T}_N^{1/2}((E_{k_1} - E_{k_1-1})\varepsilon_{k_1}(z_1))\mathbf{x}_n\right|^4\right]^{1/2}$$

$$\times \left[E\left|\sum_{k_2 \neq j} \mathbf{x}_n^*((E_{k_2} - E_{k_2-1})\varepsilon_{k_2}(z_1))\mathbf{T}_N^{1/2}\mathbf{e}_i\right|^4\right]^{1/2}$$

$$+ M \sum_{i=1}^{N} |\mathbf{x}_n^* E\mathbf{A}_j^{-1}(z_1)\mathbf{T}_N^{1/2}\mathbf{e}_i|^2 E\left|\sum_{k \neq j} \mathbf{e}_i^* \mathbf{T}_N^{1/2}((E_k - E_{k-1})\varepsilon_k(z_1))\mathbf{x}_n\right|^2$$

$$\leq M\varepsilon_N^4 + \frac{M}{N},$$

which implies

$$\frac{1}{K} \sum_{j=1}^{K} \sum_{i=1}^{N} E_j(\mathbf{H}_{nj}(z_1) - \mathbf{T}_N^{1/2}(E\mathbf{A}_j^{-1}(z))\mathbf{x}_n\mathbf{x}_n^*\mathbf{A}_j^{-1}(z)\mathbf{T}_N^{1/2})_{ii}$$
(2.2)
$$\times (E_j\mathbf{H}_{nj}(z_2))_{ii} \xrightarrow{i.p.} 0.$$

Similarly, one can also prove that

$$\frac{1}{K} \sum_{j=1}^{K} \sum_{i=1}^{N} (\mathbf{T}_N^{1/2}(E\mathbf{A}_j^{-1}(z_1))\mathbf{x}_n\mathbf{x}_n^* E_j(\mathbf{A}_j^{-1}(z_1) - E\mathbf{A}_j^{-1}(z_1))\mathbf{T}_N^{1/2})_{ii}$$

$$\times (E_j \mathbf{H}_{nj}(z_2))_{ii} \xrightarrow{i.p.} 0$$

and, therefore,

$$\frac{1}{K}\sum_{j=1}^{K}\sum_{i=1}^{N} E_j(\mathbf{H}_{nj}(z_1) - E\hat{\mathbf{H}}_{nj}(z_1))_{ii}(E_j \mathbf{H}_{nj}(z_2))_{ii} \xrightarrow{i.p.} 0.$$

Via an analogous argument,

$$\frac{1}{K}\sum_{j=1}^{K}\sum_{i=1}^{N} (E\hat{\mathbf{H}}_{nj}(z_1))_{ii} E_j(\mathbf{H}_{nj}(z_2) - E\hat{\mathbf{H}}_{nj}(z_2))_{ii} \xrightarrow{i.p.} 0.$$

Thus, for the proof of (2.1), it is sufficient to show that

(2.3) $$\sum_{i=1}^{N} (E\hat{\mathbf{H}}_{n1}(z_1))_{ii}(E\hat{\mathbf{H}}_{n1}(z_2))_{ii} \xrightarrow{i.p.} 0.$$

To this end, write

$$\mathbf{A}_1(z) - (-\hat{\mathbf{T}}_N(z)) = \sum_{k=2}^{K} \mathbf{s}_k \mathbf{s}_k^* - (-zE\underline{m}_n(z))\mathbf{T}_N,$$

where $\underline{m}_n(z)$ denotes the Stieljes transform of $\frac{N}{K}\mathbf{S}_1^* \mathbf{T}_N \mathbf{S}_1$ and $\hat{\mathbf{T}}_N(z) = zE\underline{m}_n(z) \times \mathbf{T}_N + z\mathbf{I}$. Using equality, similar to (2.2) of Silverstein [15],

(2.4) $$\underline{m}_n(z) = -\frac{1}{zK}\sum_{k=2}^{K}\beta_{1k}(z),$$

we get

$$E\mathbf{A}_1^{-1}(z) - (-\hat{\mathbf{T}}_N(z))^{-1}$$

$$= (\hat{\mathbf{T}}_N(z))^{-1} E\left[\left(\sum_{k=2}^{K}\mathbf{s}_k\mathbf{s}_k^* - (-zE\underline{m}_n(z))\mathbf{T}_N\right)\mathbf{A}_1^{-1}(z)\right]$$

(2.5)
$$= \sum_{k=2}^{K} E\left[\beta_{1k}(z)\left[(\hat{\mathbf{T}}_N(z))^{-1}\mathbf{s}_k\mathbf{s}_k^*\mathbf{A}_{1k}^{-1}(z)\right.\right.$$
$$\left.\left. - \frac{1}{K}(\hat{\mathbf{T}}_N(z))^{-1}\mathbf{T}_N E\mathbf{A}_1^{-1}(z)\right]\right].$$

It follows that

$$\mathbf{e}_i^* \mathbf{T}_N^{1/2} E\mathbf{A}_1^{-1}(z)\mathbf{x}_n - \mathbf{e}_i^* \mathbf{T}_N^{1/2}(-\hat{\mathbf{T}}_N(z))^{-1}\mathbf{x}_n$$

(2.6)
$$= (K-1)E\left[\beta_{12}(z)\left[\mathbf{s}_2^* \mathbf{A}_{12}^{-1}(z)\mathbf{x}_n \mathbf{e}_i^* \mathbf{T}_N^{1/2}(\hat{\mathbf{T}}_N(z))^{-1}\mathbf{s}_2\right.\right.$$

$$-\frac{1}{K}\mathbf{e}_i^*\mathbf{T}_N^{1/2}(\hat{\mathbf{T}}_N(z))^{-1}\mathbf{T}_N E\mathbf{A}_1^{-1}(z)\mathbf{x}_n\bigg]\bigg]$$

$$=\rho_1+\rho_2+\rho_3,$$

where

$$\rho_1 = (K-1)E[\beta_{12}(z)b_{12}(z)\xi(z)\alpha(z)],$$

$$\rho_2 = \frac{K-1}{K}E[\beta_{12}(z)\mathbf{e}_i^*\mathbf{T}_N^{1/2}(\hat{\mathbf{T}}_N(z))^{-1}\mathbf{T}_N(\mathbf{A}_{12}^{-1}(z)-\mathbf{A}_1^{-1}(z))\mathbf{x}_n]$$

and

$$\rho_3 = \frac{K-1}{K}E[\beta_{12}(z)\mathbf{e}_i^*\mathbf{T}_N^{1/2}(\hat{\mathbf{T}}_N(z))^{-1}\mathbf{T}_N(\mathbf{A}_1^{-1}(z)-E\mathbf{A}_1^{-1}(z))\mathbf{x}_n].$$

Here we also set

$$\xi(z) = \mathbf{s}_2^*\mathbf{A}_{12}^{-1}(z)\mathbf{x}_n\mathbf{e}_i^*\mathbf{T}_N^{1/2}(\hat{\mathbf{T}}_N(z))^{-1}\mathbf{s}_2$$

$$-\frac{1}{K}\mathbf{e}_i^*\mathbf{T}_N^{1/2}(\hat{\mathbf{T}}_N(z))^{-1}\mathbf{T}_N\mathbf{A}_{12}^{-1}(z)\mathbf{x}_n$$

and

$$\alpha(z) = \mathbf{s}_2^*\mathbf{A}_{12}^{-1}(z)\mathbf{s}_2 - \frac{1}{K}\operatorname{Tr}\mathbf{A}_{12}^{-1}(z), \qquad b_{12}(z) = \frac{1}{1+(1/K)\operatorname{Tr}\mathbf{A}_{12}^{-1}(z)}.$$

According to (4.2) and (4.3) in Bai, Miao and Pan [5], one can conclude that

$$\max_i|\rho_1| = O(K^{-1/2}),$$

$$\max_i|\rho_2| = \max_i\frac{K-1}{K}|E[\beta_{12}^2(z)\mathbf{e}_i^*\mathbf{T}_N^{1/2}(\hat{\mathbf{T}}_N(z))^{-1}$$

$$\times \mathbf{T}_N\mathbf{A}_{12}^{-1}(z)\mathbf{s}_2\mathbf{s}_2^*\mathbf{A}_{12}^{-1}(z)\mathbf{x}_n]|$$

$$= O(K^{-1})$$

and

$$\max_i|\rho_3| = \max_i\bigg|\frac{K-1}{K}E[\beta_{12}(z)b_{12}(z)\alpha(z)\mathbf{e}_i^*\mathbf{T}_N^{1/2}(\hat{\mathbf{T}}_N(z))^{-1}$$

$$\times \mathbf{T}_N(\mathbf{A}_1^{-1}(z)-E\mathbf{A}_1^{-1}(z))\mathbf{x}_n]\bigg|$$

$$= O(K^{-1/2}).$$

Hence,

$$\max_i|\mathbf{e}_i^*\mathbf{T}_N^{1/2}E\mathbf{A}_1^{-1}(z)\mathbf{x}_N| \to 0,$$

which, together with the Hölder inequality, guarantees (2.3). Thus, we are done.

**3. Proof of Theorem 1.1.** It is easy to show that

$$\mathbf{s}_1^*\mathbf{R}_1^m\mathbf{s}_1 - a_m \xrightarrow{i.p.} 0.$$

It follows that

(3.1) $$\mathbf{s}_1^*\mathbf{A}_{1m} - \mathbf{b}^* \xrightarrow{i.p.} 0, \qquad \mathbf{A}_{1m}^*\mathbf{R}_1\mathbf{A}_{1m} - \mathbf{B} \xrightarrow{i.p.} 0.$$

It is then observed that

$$\begin{aligned}
\sqrt{N}&(\beta_{1m} - \mathbf{b}^*\mathbf{B}^{-1}\mathbf{b}) \\
&= \sqrt{N}(\mathbf{s}_1^*\mathbf{A}_{1m} - \mathbf{b}^*)(\mathbf{A}_{1m}^*\mathbf{R}_1\mathbf{A}_{1m})^{-1}\mathbf{A}_{1m}^*\mathbf{s}_1 \\
&\quad + \sqrt{N}\mathbf{b}^*(\mathbf{A}_{1m}^*\mathbf{R}_1\mathbf{A}_{1m})^{-1}(\mathbf{A}_{1m}^*\mathbf{s}_1 - \mathbf{b}) \\
&\quad + \sqrt{N}\mathbf{b}^*((\mathbf{A}_{1m}^*\mathbf{R}_1\mathbf{A}_{1m})^{-1} - \mathbf{B}^{-1})\mathbf{b} \\
&= 2\sqrt{N}(\mathbf{s}_1^*\mathbf{A}_{1m} - \mathbf{b}^*)\mathbf{B}^{-1}\mathbf{b} \\
&\quad - \sqrt{N}\mathbf{b}^*\mathbf{B}^{-1}(\mathbf{A}_{1m}^*\mathbf{R}_1\mathbf{A}_{1m} - \mathbf{B})\mathbf{B}^{-1}\mathbf{b} + o_p(1),
\end{aligned}$$ (3.2)

where we use (3.1), (3.6) below and an identity

$$\mathbf{B}_1^{-1} - \mathbf{B}_2^{-1} = -\mathbf{B}_1^{-1}(\mathbf{B}_1 - \mathbf{B}_2)\mathbf{B}_2^{-1},$$

which holds for any invertible matrices $\mathbf{B}_1$ and $\mathbf{B}_2$. Furthermore, let

$$\mathbf{b}^*\mathbf{B}^{-1} = (d_1, \ldots, d_m),$$

then (3.2) is now equal to

(3.3) $$2\sqrt{N}\sum_{i=1}^m d_i(\mathbf{s}_1^*\mathbf{R}_1^{i-1}\mathbf{s}_1 - a_{i-1}) - \sqrt{N}\sum_{i,j=1}^m d_i d_j(\mathbf{s}_1^*\mathbf{R}_1^{i+j-1}\mathbf{s}_1 - a_{i+j-1}).$$

By the result (1) of Theorem 1.1 of Bai and Silverstein [3], it is easily seen that

$$\sqrt{N}\left(\frac{1}{N}\operatorname{Tr}\mathbf{R}_1^i - a_i\right) \xrightarrow{i.p.} 0.$$

To derive a central limit theorem for (3.3), it then suffices to develop a multivariate one for $\{\sqrt{N}(\mathbf{s}_1^*\mathbf{R}_1^i\mathbf{s}_1 - \frac{1}{N}\operatorname{Tr}\mathbf{R}_1^i), i = 0, \ldots, 2m-1\}$.

Set $\mathbf{H}_1 = \mathbf{S}_1\mathbf{S}_1^*$ and $h_m = \int x^m \, dF^{c_N}(cx)$. Note that

(3.4) $$\sqrt{N}\left(\mathbf{s}_1^*\mathbf{R}_1^i\mathbf{s}_1 - \frac{1}{N}\operatorname{Tr}\mathbf{R}^i\right) = \sum_{u=0}^i \binom{i}{u}(\sigma^2)^{i-u}\sqrt{N}\left(\mathbf{s}_1^*\mathbf{H}_1^u\mathbf{s}_1 - \frac{1}{N}\operatorname{Tr}\mathbf{H}_1^u\right).$$

Let $\|\mathbf{s}_1\|^2 = \sum_{i=1}^{N} |v_{i1}|^2/N$. Write

$$\sqrt{N}\left(\mathbf{s}_1^*\mathbf{H}_1^u\mathbf{s}_1 - \frac{1}{N}\operatorname{Tr}\mathbf{H}_1^u\right)$$
$$= \sqrt{N}\|\mathbf{s}_1\|^2\left(\frac{\mathbf{s}_1^*\mathbf{H}_1^u\mathbf{s}_1}{\|\mathbf{s}_1\|^2} - \frac{1}{N}\operatorname{Tr}\mathbf{H}_1^u\right) + \sqrt{N}\frac{1}{N}\operatorname{Tr}\mathbf{H}_1^u(\|\mathbf{s}_1\|^2 - 1).$$

It is easy to check that

$$\max_i \left|\frac{v_{i1}/\sqrt{N}}{\|\mathbf{s}_1\|}\right| \xrightarrow{i.p.} 0.$$

Therefore, given $\mathbf{s}_1$, it follows from Theorem 1.3 that

$$\left(\sqrt{N}\left(\frac{\mathbf{s}_1^*\mathbf{H}_1\mathbf{s}_1}{\|\mathbf{s}_1\|^2} - \frac{1}{N}\operatorname{Tr}\mathbf{H}_1\right), \ldots, \sqrt{N}\left(\frac{\mathbf{s}_1^*\mathbf{H}_1^u\mathbf{s}_1}{\|\mathbf{s}_1\|^2} - \frac{1}{N}\operatorname{Tr}\mathbf{H}_1^u\right)\right)$$
$$\xrightarrow{D} \sqrt{2}\left(\int_{(1-\sqrt{c})^2}^{(1+\sqrt{c})^2} \frac{x}{c} dW_x^c, \ldots, \int_{(1-\sqrt{c})^2}^{(1+\sqrt{c})^2} \frac{x^u}{c^u} dW_x^c\right)$$

(regarding the formula, one may refer to Bai, Miao and Pan [5] or Silverstein [13, 14]). However, it is evident that

(3.5)  $$\sqrt{N}(\|\mathbf{s}_1\|^2 - 1) \xrightarrow{D} X,$$

where $X \sim N(0, E|v_{11}|^4 - 1)$. Consequently, by the independence of $\mathbf{s}_1$ and $\mathbf{H}_1$,

$$\left(\sqrt{N}(\mathbf{s}_1^*\mathbf{s}_1 - 1), \ldots, \sqrt{N}\left(\mathbf{s}_1^*\mathbf{H}_1^{2m-1}\mathbf{s}_1 - \frac{1}{N}\operatorname{Tr}\mathbf{H}_1^{2m-1}\right)\right) \xrightarrow{D} (\xi_0, \ldots, \xi_{2m-1}),$$

where $\xi_i = h_i X + \frac{\sqrt{2}}{c^i}\int_{(1-\sqrt{c})^2}^{(1+\sqrt{c})^2} x^i dW_x^c, i = 1, \ldots, 2m-1$, and $\xi_0 = X$. Then

(3.6)
$$\left(\sqrt{N}(\mathbf{s}_1^*\mathbf{s}_1 - 1), \ldots, \sqrt{N}\left(\mathbf{s}_1^*\mathbf{R}_1^{2m-1}\mathbf{s}_1 - \frac{1}{N}\operatorname{Tr}\mathbf{R}_1^{2m-1}\right)\right)$$
$$\xrightarrow{D} (\zeta_0, \ldots, \zeta_{2m-1}),$$

where $\zeta_i = \sum_{u=0}^{i} \binom{i}{u}(\sigma^2)^{i-u}\xi_u$.

It follows that

$$\sqrt{N}(\beta_{1m} - \mathbf{b}^*\mathbf{B}^{-1}\mathbf{b}) \xrightarrow{D} 2\sum_{i=1}^{m} d_i\zeta_i - \sum_{i,j=1}^{m} d_i d_j \zeta_{i+j-1}.$$

Thus, we are done.

**4. Proof of Theorem 1.4.** By the argument of Bai and Silverstein [3], it suffices to find the limits of the following sums:

$$(4.1) \quad \frac{1}{K^2} \sum_{j=1}^{K} \sum_{i=1}^{N} E_j(\mathbf{T}_N^{1/2} \mathbf{A}_j^{-1}(z_1) \mathbf{T}_N^{1/2})_{ii} E_j(\mathbf{T}_N^{1/2} \mathbf{A}_j^{-1}(z_2) \mathbf{T}_N^{1/2})_{ii}$$

and

$$(4.2) \quad \frac{1}{K} \sum_{i=1}^{N} E[(\mathbf{T}_N^{1/2} \mathbf{A}_j^{-1}(z) \mathbf{T}_N^{1/2})_{ii} (\mathbf{T}_N^{1/2} \mathbf{A}_j^{-1}(z)(\hat{\mathbf{T}}_N(z))^{-1} \mathbf{T}_N^{1/2})_{ii}]$$

(see (2.7) and (4.10) in Bai and Silverstein [3]).

Similar to (2.2), it can be verified that

$$\frac{1}{K^2} \sum_{j=1}^{K} \sum_{i=1}^{N} E_j(\mathbf{T}_N^{1/2} \mathbf{A}_j^{-1}(z_1) \mathbf{T}_N^{1/2} - E(\mathbf{T}_N^{1/2} \mathbf{A}_j^{-1}(z_1) \mathbf{T}_N^{1/2}))_{ii}$$

$$\times E_j(\mathbf{T}_N^{1/2} \mathbf{A}_j^{-1}(z_2) \mathbf{T}_N^{1/2})_{ii} = O_p(N^{-1/2}).$$

Consequently, analogous to Theorem 1.3, it remains to find the limit of

$$(4.3) \quad \frac{1}{K} \sum_{i=1}^{N} E(\mathbf{T}_N^{1/2} \mathbf{A}_1^{-1}(z_1) \mathbf{T}_N^{1/2})_{ii} E(\mathbf{T}_N^{1/2} \mathbf{A}_1^{-1}(z_2) \mathbf{T}_N^{1/2})_{ii}.$$

Define

$$\gamma(z) = \mathbf{s}_k^* \mathbf{A}_{1k}^{-1}(z) \mathbf{T}_N^{1/2} \mathbf{e}_i \mathbf{e}_i^* \mathbf{T}_N^{1/2} (\hat{\mathbf{T}}_N(z))^{-1} \mathbf{s}_k$$

$$- \frac{1}{K} \mathbf{e}_i^* \mathbf{T}_N^{1/2} (\hat{\mathbf{T}}_N(z))^{-1} \mathbf{T}_N \mathbf{A}_{1k}^{-1}(z) \mathbf{T}_N^{1/2} \mathbf{e}_i.$$

From (2.5), we have

$$(4.4) \quad \begin{aligned} E(\mathbf{T}_N^{1/2} \mathbf{A}_1^{-1}(z) \mathbf{T}_N^{1/2})_{ii} &- \mathbf{e}_i^* \mathbf{T}_N^{1/2} (-\hat{\mathbf{T}}_N(z))^{-1} \mathbf{T}_N^{1/2} \mathbf{e}_i \\ &= \sum_{k=2}^{K} E\bigg[ \beta_{1k}(z) \mathbf{e}_i^* \mathbf{T}_N^{1/2} (\hat{\mathbf{T}}_N(z))^{-1} \mathbf{s}_k \mathbf{s}_k^* \mathbf{A}_{1k}^{-1}(z) \mathbf{T}_N^{1/2} \mathbf{e}_i \\ &\quad - \beta_{1k}(z) \mathbf{e}_i^* \mathbf{T}_N^{1/2} \frac{1}{K} (\hat{\mathbf{T}}_N(z_1))^{-1} \mathbf{T}_N E \mathbf{A}_1^{-1}(z) \mathbf{T}_N^{1/2} \mathbf{e}_i \bigg] \\ &= \tau_1(z) + \tau_2(z) + \tau_3(z), \end{aligned}$$

where

$$\tau_1(z) = (K-1) E[\beta_{12}(z) b_{12}(z) \gamma(z) \alpha(z)],$$

$$\tau_2(z) = \frac{K-1}{K} E[\beta_{12}(z) \mathbf{e}_i^* \mathbf{T}_N^{1/2} (\hat{\mathbf{T}}_N(z))^{-1} \mathbf{T}_N (\mathbf{A}_{12}^{-1}(z) - \mathbf{A}_1^{-1}(z)) \mathbf{T}_N^{1/2} \mathbf{e}_i]$$

and

$$\tau_3(z) = \frac{K-1}{K} E[\beta_{12}(z) \mathbf{e}_i^* \mathbf{T}_N^{1/2} (\hat{\mathbf{T}}_N(z))^{-1} \mathbf{T}_N (\mathbf{A}_1^{-1}(z) - E\mathbf{A}_1^{-1}(z)) \mathbf{T}_N^{1/2} \mathbf{e}_i].$$

Therefore, it follows from (4.4) that

$$\frac{1}{K} \sum_{i=1}^{N} E(\mathbf{T}_N^{1/2} \mathbf{A}_1^{-1}(z_1) \mathbf{T}_N^{1/2})_{ii} E(\mathbf{T}_N^{1/2} \mathbf{A}_1^{-1}(z_2) \mathbf{T}_N^{1/2})_{ii}$$

$$= \frac{1}{K} \sum_{i=1}^{N} \mathbf{e}_i^* \mathbf{T}_N^{1/2} (\hat{\mathbf{T}}_N(z_1))^{-1} \mathbf{T}_N^{1/2} \mathbf{e}_i \mathbf{e}_i^* \mathbf{T}_N^{1/2} (\hat{\mathbf{T}}_N(z_2))^{-1} \mathbf{T}_N^{1/2} \mathbf{e}_i + O\left(\frac{1}{\sqrt{K}}\right),$$

where the estimate can be obtained as in Theorem 1.3.

Regarding (4.2), due to similar reason, one need only seek the limit of

$$\frac{1}{K} \sum_{i=1}^{N} E[(\mathbf{T}_N^{1/2} \mathbf{A}_j^{-1}(z) \mathbf{T}_N^{1/2})_{ii}] E[(\mathbf{T}_N^{1/2} \mathbf{A}_j^{-1}(z) (\hat{\mathbf{T}}_N(z))^{-1} \mathbf{T}_N^{1/2})_{ii}].$$

However, as in (4.4), one can conclude that

$$\frac{1}{K} \sum_{i=1}^{N} E[(\mathbf{T}_N^{1/2} \mathbf{A}_j^{-1}(z) \mathbf{T}_N^{1/2})_{ii}] E[(\mathbf{T}_N^{1/2} \mathbf{A}_j^{-1}(z) (\hat{\mathbf{T}}_N(z))^{-1} \mathbf{T}_N^{1/2})_{ii}]$$

$$= \frac{1}{K} \sum_{i=1}^{N} \mathbf{e}_i^* \mathbf{T}_N^{1/2} (\hat{\mathbf{T}}_N(z))^{-1} \mathbf{T}_N^{1/2} \mathbf{e}_i \mathbf{e}_i^* \mathbf{T}_N^{1/2} (\hat{\mathbf{T}}_N(z))^{-2} \mathbf{T}_N^{1/2} \mathbf{e}_i + O\left(\frac{1}{\sqrt{K}}\right).$$

For later purpose, we now derive (1.23) and (1.24). Note that when $\mathbf{T}_N = \mathbf{I}$, for $z \in \mathbf{C}^+$,

(4.5) $$z = -\frac{1}{\underline{m}(z)} + \frac{c}{1 + \underline{m}(z)}$$

and

(4.6) $$\frac{d}{dz}\underline{m}(z) = \frac{\underline{m}^2(z)}{1 - c\underline{m}^2(z)/(1 + \underline{m}(z))^2}.$$

Then for $g(x) = x^r$,

$$\frac{1}{2\pi i} \int g(z) \frac{c\underline{m}^3(z) h_2(z)}{1 - c \int \underline{m}^2(z) t^2 \, dH(t)/(1 + t\underline{m}(z))^2} \, dz$$

$$= \frac{c}{2\pi i} \int \frac{(-1/\underline{m}(z) + c/(1 + \underline{m}(z)))^r}{(\underline{m}(z) + 1)^3} \underline{m}(z) \, d\underline{m}(z)$$

$$= \frac{c}{2\pi i} \int \frac{(-1/\underline{m}(z) + c/(1 + \underline{m}(z)))^r}{(\underline{m}(z) + 1)^2} \, d\underline{m}(z)$$

$$-\frac{c}{2\pi i}\int\frac{(-1/\underline{m}(z)+c/(1+\underline{m}(z)))^r}{(\underline{m}(z)+1)^3}d\underline{m}(z)$$

$$\stackrel{\triangle}{=}\nu_1-\nu_2.$$

For $\nu_1$, we have

$$\nu_1=c^r\frac{c}{2\pi i}\int\frac{((1-c)/c+1/(1+\underline{m}(z)))^r}{(\underline{m}(z)+1)^2}(1-(1+\underline{m}(z)))^{-r}d\underline{m}(z)$$

$$=\frac{c^{1+r}}{2\pi i}\int\sum_{j=0}^{r}\binom{r}{j}\left(\frac{1-c}{c}\right)^j\frac{1}{(1+\underline{m}(z))^{r-j+2}}$$

$$\times\sum_{k=0}^{\infty}\binom{r+k-1}{k}(1+\underline{m}(z))^k\,d\underline{m}(z)$$

$$=c^{1+r}\sum_{j=0}^{r}\binom{r}{j}\left(\frac{1-c}{c}\right)^j\binom{2r-j}{r-1}.$$

Similarly,

$$\nu_2=c^{1+r}\sum_{j=0}^{r}\binom{r}{j}\left(\frac{1-c}{c}\right)^j\binom{2r+1-j}{r-1}.$$

For (1.24), we have

$$\int z_1^{r_1}\frac{d}{dz_1}\left[\frac{\underline{m}(z_1)}{1+\underline{m}(z_1)}\right]dz_1=\int\frac{(-1/\underline{m}(z_1)+c/(1+\underline{m}(z_1)))^r}{(\underline{m}(z_1)+1)^2}d\underline{m}(z_1)$$

$$=2\pi ic^{r_1}\sum_{j=0}^{r_1}\binom{r_1}{j}\left(\frac{1-c}{c}\right)^j\binom{2r_1-j}{r_1-1}.$$

Therefore,

$$-\frac{c}{4\pi^2}\int\int g_1(z_1)g_2(z_2)\frac{d^2}{dz_1\,dz_2}[\underline{m}(z_1)\underline{m}(z_2)h_1(z_1,z_2)]\,dz_1\,dz_2$$

$$=c^{r_1+r_2+1}\sum_{j_1=0}^{r_1}\sum_{j_2=0}^{r_2}\binom{r_1}{j_1}\binom{r_2}{j_2}\left(\frac{1-c}{c}\right)^{j_1+j_2}\binom{2r_1-j_1}{r_1-1}\binom{2r_2-j_2}{r_2-1}.$$

**5. Proof of Theorem 1.2.** Since the truncation process is tedious, it is deferred to the Appendix. It may then be assumed that the underlying r.v.'s satisfy

$$Ev_{11}=0,\qquad E|v_{11}|^2=1,\qquad |v_{11}|\le\varepsilon_N\sqrt{N},$$

where $\varepsilon_N$ is a positive sequence converging to zero.

Define $\check{s}_k = \mathbf{s}_k^* \mathbf{R}_k \mathbf{s}_k - a_1$. Expand $(\mathbf{s}_k^* \mathbf{R}_k \mathbf{s}_k)^{-1}$ a little bit as follows:

$$\frac{1}{\mathbf{s}_k^* \mathbf{R}_k \mathbf{s}_k} = \frac{1}{a_1} - \frac{\check{s}_k}{a_1 \mathbf{s}_k^* \mathbf{R}_k \mathbf{s}_k}$$

(5.1)

$$= \frac{1}{a_1} - \frac{\check{s}_k}{a_1^2} + \frac{(\check{s}_k)^2}{a_1^2 \mathbf{s}_k^* \mathbf{R}_k \mathbf{s}_k}.$$

It follows that

(5.2) $$\sum_{k=1}^{K} \left( \beta_k - \frac{p_k}{\alpha_1} \right) = G_1 + G_2 + G_3 + G_4,$$

where

$$G_1 = \frac{1}{a_1} \sum_{k=1}^{K} p_k((\mathbf{s}_k^* \mathbf{s}_k)^2 - 1), \qquad G_2 = -\frac{1}{a_1^2} \sum_{k=1}^{K} p_k(\mathbf{s}_k^* \mathbf{s}_k)^2 (\check{s}_k)$$

and

$$G_3 = \frac{1}{a_1^3} \sum_{k=1}^{K} p_k(\mathbf{s}_k^* \mathbf{s}_k)^2 (\check{s}_k)^2, \qquad G_4 = -\frac{1}{a_1^3} \sum_{k=1}^{K} \frac{p_k(\mathbf{s}_k^* \mathbf{s}_k)^2 (\check{s}_k)^3}{\mathbf{s}_k^* \mathbf{R}_k \mathbf{s}_k}.$$

We will analyze $G_1, G_2, G_3, G_4$ one by one and, as will be seen, the contribution from the term $G_4$ is negligible.

First consider the term $G_4$. Since $\mathbf{s}_k^* \mathbf{R}_k \mathbf{s}_k \geq \sigma^2 \mathbf{s}_k^* \mathbf{s}_k$, we have

$$|G_4| \leq M(G_{41} + \cdots + G_{43}),$$

where

$$G_{41} = \sum_{k=1}^{K} p_k \left| \mathbf{s}_k^* \mathbf{s}_k \left( \mathbf{s}_k^* \mathbf{R}_k \mathbf{s}_k - \frac{1}{N} \operatorname{Tr} \mathbf{R}_k \right)^3 \right|$$

and

$$G_{42} = \sum_{k=1}^{K} p_k \left| \mathbf{s}_k^* \mathbf{s}_k \left( \frac{1}{N} \operatorname{Tr} \mathbf{R}_k - \operatorname{Tr} \mathbf{R} \right)^3 \right|, \qquad G_{43} = \sum_{k=1}^{K} p_k \left| \mathbf{s}_k^* \mathbf{s}_k \left( \frac{1}{N} \operatorname{Tr} \mathbf{R} - a_1 \right)^3 \right|.$$

By the Hölder inequality,

$$EG_{41} \leq M \sum_{k=1}^{K} p_k (E(\mathbf{s}_k^* \mathbf{s}_k - 1)^2)^{1/2} \left( E \left( \mathbf{s}_k^* \mathbf{R}_k \mathbf{s}_k - \frac{1}{N} \operatorname{Tr} \mathbf{R}_k \right)^6 \right)^{1/2}$$

$$+ M \sum_{k=1}^{K} p_k E \left| \mathbf{s}_k^* \mathbf{R}_k \mathbf{s}_k - \frac{1}{N} \operatorname{Tr} \mathbf{R}_k \right|^3$$

$$= o(1).$$

Indeed, it is easy to verify that

$$E(\mathbf{s}_k^*\mathbf{s}_k - 1)^2 = \frac{1}{N}(E|v_{11}|^4 - 1) \tag{5.3}$$

and that

$$E\left(\mathbf{s}_k^*\mathbf{R}_k\mathbf{s}_k - \frac{1}{N}\operatorname{Tr}\mathbf{R}_k\right)^p \leq \frac{M}{N^p}(E|v_{11}|^4 E(\operatorname{Tr}\mathbf{R}_k^2)^{p/2} + Ev_{11}^{2p} E\operatorname{Tr}\mathbf{R}_k^p)$$

$$\leq \frac{M}{N^{p/2}} + \frac{M\varepsilon_N^{2p-4}}{N^2}, \tag{5.4}$$

where the constant $M$ is independent of $k$. Here we use the fact $\mathbf{R}_k \leq M\mathbf{S}_k\mathbf{S}_k^* + \sigma^2\mathbf{I}$.

Furthermore, it is direct to prove

$$\frac{1}{N^2}\sum_{k=1}^K p_k|\mathbf{s}_k^*\mathbf{s}_k| \xrightarrow{i.p.} 0.$$

This, together with Theorem 1 of Bai and Silverstein [3], leads to

$$G_{43} \xrightarrow{i.p.} 0.$$

In addition, it is also easy to verify that

$$EG_{42} = \frac{1}{N^3}\sum_{k=1}^K p_k^4 E(\mathbf{s}_k^*\mathbf{s}_k)^4 = O\left(\frac{1}{N^2}\right).$$

Combining the above argument, one can claim that the contribution from $G_4$ can be ignored.

Analyze the term $G_1$ second. Write

$$\sum_{k=1}^K p_k(\mathbf{s}_k^*\mathbf{s}_k)^2 = \sum_{k=1}^K p_k(\mathbf{s}_k^*\mathbf{s}_k - 1)^2 + 2\sum_{k=1}^K p_k\mathbf{s}_k^*\mathbf{s}_k - \sum_{k=1}^K p_k$$

$$= \sum_{k=1}^K p_k(\mathbf{s}_k^*\mathbf{s}_k - 1)^2 + 2\operatorname{Tr}\mathbf{C} - \sum_{k=1}^K p_k. \tag{5.5}$$

Moreover,

$$E\left(\sum_{k=1}^K [p_k(\mathbf{s}_k^*\mathbf{s}_k - 1)^2 - E(\mathbf{s}_k^*\mathbf{s}_k - 1)^2]\right)^2$$

$$= \sum_{k=1}^K p_k^2 E((\mathbf{s}_k^*\mathbf{s}_k - 1)^2 - E(\mathbf{s}_k^*\mathbf{s}_k - 1)^2)^2 = o(1),$$

using

(5.6) $$E(\mathbf{s}_k^*\mathbf{s}_k - 1)^4 = o\left(\frac{1}{N}\right).$$

So

$$\sum_{k=1}^{K} p_k(\mathbf{s}_k^*\mathbf{s}_k - 1)^2 \xrightarrow{i.p.} \frac{1}{c}(E|v_{11}|^4 - 1)\int x\,dH(x),$$

and then

(5.7) $$G_1 = \frac{1}{a_1}\left(\frac{(E|v_{11}|^4 - 1)\int x\,dH(x)}{c} + 2\operatorname{Tr}\mathbf{C} - 2\sum_{k=1}^{K} p_k\right) + o_p(1).$$

Third, for the term $G_2$, similar to $G_1$,

(5.8) $$-a_1^2 G_2 = \sum_{k=1}^{K} p_k(\check{s}_k)(\mathbf{s}_k^*\mathbf{s}_k - 1)^2$$

(5.9) $$+ 2\sum_{k=1}^{K} p_k(\check{s}_k)(\mathbf{s}_k^*\mathbf{s}_k - 1) + \sum_{k=1}^{K} p_k(\check{s}_k).$$

For the sum in (5.8), we have

$$E\left|\sum_{k=1}^{K} p_k(\check{s}_k)(\mathbf{s}_k^*\mathbf{s}_k - 1)^2\right| \leq M \sum_{k=1}^{K} (E(\check{s}_k)^2)^{1/2}(E(\mathbf{s}_k^*\mathbf{s}_k - 1)^4)^{1/2}$$
$$= o(1),$$

where we use (5.6) and

$$E(\check{s}_k)^2 \leq E\left(\mathbf{s}_k^*\mathbf{R}_k\mathbf{s}_k - \frac{1}{N}\operatorname{Tr}\mathbf{R}_k\right)^2 + E\left(\frac{1}{N}\operatorname{Tr}\mathbf{R}_k - a_1\right)^2 = O\left(\frac{1}{N}\right),$$

which is accomplished by (5.4) and Theorem 1 of Bai and Silverstein [3]. Similarly to (5.5), we deduce that

(5.10) $$\sum_{k=1}^{K} p_k^2(\mathbf{s}_k^*\mathbf{s}_k)^2 = \frac{1}{c}(E|v_{11}|^4 - 1)\int x^2\,dH(x)$$
$$+ 2\operatorname{Tr}\mathbf{S}\mathbf{P}^2\mathbf{S}^* - \sum_{k=1}^{K} p_k^2 + o_p(1).$$

Applying $\mathbf{C} - p_k\mathbf{s}_k\mathbf{s}_k^* = \mathbf{C}_k$, the second sum of (5.9) is then equal to

$$\sigma^2\operatorname{Tr}\mathbf{C} + \operatorname{Tr}\mathbf{C}^2 - a_1\sum_{k=1}^{K} p_k - \sum_{k=1}^{K} p_k^2(\mathbf{s}_k^*\mathbf{s}_k)^2$$

$$\text{(5.11)} \quad = \sigma^2 \operatorname{Tr} \mathbf{C} + \operatorname{Tr} \mathbf{C}^2 - a_1 \sum_{k=1}^{K} p_k - \frac{1}{c}(E|v_{11}|^4 - 1) \int x^2 \, dH(x)$$

$$- 2 \operatorname{Tr} \mathbf{S}\mathbf{P}^2\mathbf{S}^* + \sum_{k=1}^{K} p_k^2 + o_p(1).$$

With regard to the first sum of (5.9), its variance will be proved to converge to zero. Now let us provide more details to the reader:

$$\text{(5.12)} \quad \operatorname{Var}\left(\sum_{k=1}^{K} p_k(\check{s}_k)(\mathbf{s}_k^*\mathbf{s}_k - 1)\right) = G_{21} + G_{22},$$

where

$$G_{21} = \sum_{k=1}^{K} p_k^2 E[(\check{s}_k)(\mathbf{s}_k^*\mathbf{s}_k - 1) - E(\check{s}_k)(\mathbf{s}_k^*\mathbf{s}_k - 1)]^2$$

and

$$G_{22} = \sum_{k_1 \neq k_2}^{K} p_{k_1} p_{k_2} E[((\check{s}_{k_1})(\mathbf{s}_{k_1}^*\mathbf{s}_{k_1} - 1) - E(\check{s}_{k_1})(\mathbf{s}_{k_1}^*\mathbf{s}_{k_1} - 1))$$

$$\times ((\check{s}_{k_2})(\mathbf{s}_{k_2}^*\mathbf{s}_{k_2} - 1) - E(\check{s}_{k_2})(\mathbf{s}_{k_2}^*\mathbf{s}_{k_2} - 1))].$$

Evidently,

$$G_{21} \leq \sum_{k=1}^{K} M E[(\check{s}_k)(\mathbf{s}_k^*\mathbf{s}_k - 1)]^2$$

$$\leq M \sum_{k=1}^{K} E\left[\left(\mathbf{s}_k^*\mathbf{R}_k\mathbf{s}_k - \frac{1}{N}\operatorname{Tr}\mathbf{R}_k\right)(\mathbf{s}_k^*\mathbf{s}_k - 1)\right]^2$$

$$+ M \sum_{k=1}^{K} E\left[\left(\frac{1}{N}\operatorname{Tr}\mathbf{R}_k - a_1\right)(\mathbf{s}_k^*\mathbf{s}_k - 1)\right]^2$$

$$\text{(5.13)} \quad \leq M \sum_{k=1}^{K} \left[E\left(\mathbf{s}_k^*\mathbf{R}_k\mathbf{s}_k - \frac{1}{N}\operatorname{Tr}\mathbf{R}_k\right)^4\right]^{1/2} [E(\mathbf{s}_k^*\mathbf{s}_k - 1)^4]^{1/2}$$

$$+ M \sum_{k=1}^{K} E\left(\frac{1}{N}\operatorname{Tr}\mathbf{R}_k - a_1\right)^2 E(\mathbf{s}_k^*\mathbf{s}_k - 1)^2$$

$$= o(1).$$

Let $\mathbf{S}_{k_1 k_2}$ denote the matrix obtained from $\mathbf{S}_{k_1}$ by deleting the $k_2$th column and, furthermore, $\mathbf{R}_{k_1 k_2}$ and $\mathbf{C}_{k_1 k_2}$ have the same meaning. Split $\mathbf{R}_{k_1} =$ $\mathbf{R}_{k_1k_2} + p_{k_2}\mathbf{s}_{k_2}\mathbf{s}_{k_2}^*$ and $\mathbf{R}_{k_2} = \mathbf{R}_{k_1k_2} + p_{k_1}\mathbf{s}_{k_1}\mathbf{s}_{k_1}^*$. Also, for convenience, set

$$\alpha_{k_j} = \mathbf{s}_{k_j}^* \mathbf{R}_{k_1k_2} \mathbf{s}_{k_j} - a_1, \qquad \gamma_j = \mathbf{s}_{k_j}^* \mathbf{R}_{k_1k_2} \mathbf{s}_{k_j} - \frac{1}{N} \operatorname{Tr} \mathbf{R}_{k_1k_2}$$

and

$$\Upsilon_{k_j} = \mathbf{s}_{k_j}^* \mathbf{s}_{k_j} - 1, \qquad j = 1, 2.$$

$G_{22}$ is then decomposed as

$$G_{22} = G_{221} + \cdots + G_{224},$$

where

$$G_{221} = \sum_{k_1 \neq k_2}^{K} p_{k_1} p_{k_2} \operatorname{Cov}(\alpha_{k_1} \Upsilon_{k_1}, \alpha_{k_2} \Upsilon_{k_2}),$$

$$G_{222} = \sum_{k_1 \neq k_2}^{K} p_{k_1}^2 p_{k_2} \operatorname{Cov}(\alpha_{k_1} \Upsilon_{k_1}, |\mathbf{s}_{k_1}^* \mathbf{s}_{k_2}|^2 \Upsilon_{k_2}),$$

$$G_{223} = \sum_{k_1 \neq k_2}^{K} p_{k_1} p_{k_2}^2 \operatorname{Cov}(|\mathbf{s}_{k_1}^* \mathbf{s}_{k_2}|^2 \Upsilon_{k_1}, \alpha_{k_2} \Upsilon_{k_2})$$

and

$$G_{224} = \sum_{k_1 \neq k_2}^{K} p_{k_1}^2 p_{k_2}^2 \operatorname{Cov}(|\mathbf{s}_{k_1}^* \mathbf{s}_{k_2}|^2 \Upsilon_{k_1}, |\mathbf{s}_{k_1}^* \mathbf{s}_{k_2}|^2 \Upsilon_{k_2}).$$

The basic idea behind this decomposition is to produce some independent terms when $\mathbf{R}_{k_1k_2}$ is given, which is very important when estimating the order of some terms.

It is easy to check that

$$(5.14) \qquad E\left(\mathbf{s}_k^* \mathbf{R}_k \mathbf{s}_k - \frac{1}{N} \operatorname{Tr} \mathbf{R}_k\right)(\mathbf{s}_k^* \mathbf{s}_k - 1) = \frac{E|v_{11}|^4 - 1}{N^2} E \operatorname{Tr} \mathbf{R}_k,$$

and that

$$(5.15) \qquad E\left|\mathbf{s}_1^* \mathbf{D} \mathbf{s}_1 - \frac{1}{N} \operatorname{Tr} \mathbf{D}\right|^2 = \frac{1}{N^2}(E|v_{11}|^4 - 2)\sum_{i=1}^{N}[(\mathbf{D})_{ii}]^2 + \frac{1}{N^2} \operatorname{Tr} \mathbf{D}\mathbf{D}^*,$$

where $\mathbf{D}$ is any constant Hermite matrix.

This gives that $G_{221}$ is equal to

$$\sum_{k_1 \neq k_2}^{K} p_{k_1} p_{k_2} E[E(\widehat{\alpha_{k_1} \Upsilon_{k_1}} \mid \mathbf{R}_{k_1k_2}) E(\widehat{\alpha_{k_2} \Upsilon_{k_2}} \mid \mathbf{R}_{k_1k_2})]$$

$$= \sum_{k_1 \neq k_2}^{K} p_{k_1} p_{k_2} E[E(\widehat{\gamma_{k_1} \Upsilon_{k_1}} \mid \mathbf{R}_{k_1 k_2}) E(\widehat{\gamma_{k_2} \Upsilon_{k_2}} \mid \mathbf{R}_{k_1 k_2})]$$

$$= \frac{Ev_{11}^4 - 1}{N} \sum_{k_1 \neq k_2}^{K} p_{k_1} p_{k_2} E\left(\frac{1}{N} \operatorname{Tr} \mathbf{R}_{k_1 k_2} - E\frac{1}{N} \operatorname{Tr} \mathbf{R}_{k_1 k_2}\right)^2 = O\left(\frac{1}{N}\right),$$

where $\widehat{\alpha_k \Upsilon_k} = \alpha_k \Upsilon_k - E\alpha_k \Upsilon_k$, $\widehat{\gamma_k \Upsilon_k} = \gamma_k \Upsilon_k - E\gamma_k \Upsilon_k$, and we use the independence of $\mathbf{s}_{k_1}$ and $\mathbf{s}_{k_2}$, and

$$E\left(\frac{1}{N} \operatorname{Tr} \mathbf{R}_{k_1 k_2} - E\frac{1}{N} \operatorname{Tr} \mathbf{R}_{k_1 k_2}\right)^2 = \frac{1}{N^2} E\left|\sum_{j \neq k_1, k_2} p_j(\mathbf{s}_j^* \mathbf{s}_j - 1)\right|^2 \leq \frac{M}{N^2},$$

(5.16)
where $M$ is independent of $k_1, k_2$.

After some simple computations, we get

(5.17)
$$E|\mathbf{s}_{k_1}^* \mathbf{s}_{k_2}|^2 \Upsilon_{k_2} = \frac{E|v_{11}|^4 - 1}{N^2},$$

$$E(|\mathbf{s}_{k_1}^* \mathbf{s}_{k_2}|^2 \Upsilon_{k_2} \mid \mathbf{s}_{k_1}) = \frac{E|v_{11}|^4 - 1}{N^2} \mathbf{s}_{k_1}^* \mathbf{s}_{k_1},$$

and so

$$G_{222} = \sum_{k_1 \neq k_2}^{K} p_{k_1}^2 p_{k_2} E((\widehat{\alpha_{k_1} \Upsilon_{k_1}}) E[|\mathbf{s}_{k_1}^* \mathbf{s}_{k_2}|^2 \Upsilon_{k_2} - E(|\mathbf{s}_{k_1}^* \mathbf{s}_{k_2}|^2 \Upsilon_{k_2}) \mid \mathbf{s}_{k_1}])$$

$$= \frac{E|v_{11}|^4 - 1}{N^2} \sum_{k_1 \neq k_2}^{K} p_{k_1}^2 p_{k_2} (E[(\widehat{\alpha_{k_1} \Upsilon_{k_1}}) \mathbf{s}_{k_1}^* \mathbf{s}_{k_1}] + E(\widehat{\alpha_{k_1} \Upsilon_{k_1}}))$$

(5.18)
$$= \frac{E|v_{11}|^4 - 1}{N^2} \sum_{k_1 \neq k_2}^{K} p_{k_1}^2 p_{k_2} E[(\widehat{\alpha_{k_1} \Upsilon_{k_1}}) \Upsilon_{k_1}]$$

$$\leq \frac{M}{N^2} \sum_{k_1 \neq k_2}^{K} (E\alpha_{k_1}^2)^{1/2} (E\Upsilon_{k_1}^4)^{1/2}$$

$$= O\left(\frac{1}{N}\right).$$

Similarly, one can conclude that

(5.19) $$G_{223} \to 0.$$

Write

$$G_{224} = \sum_{k_1 \neq k_2}^{K} p_{k_1}^2 p_{k_2}^2 E[|\mathbf{s}_{k_1}^* \mathbf{s}_{k_2}|^4 \Upsilon_{k_1} \Upsilon_{k_2}] - \sum_{k_1 \neq k_2}^{K} p_{k_1}^2 p_{k_2}^2 [E(|\mathbf{s}_{k_1}^* \mathbf{s}_{k_2}|^2 \Upsilon_{k_1})]^2.$$

The second sum converges to zero because of (5.17). For its first sum we have

$$\sum_{k_1 \neq k_2}^{K} p_{k_1}^2 p_{k_2}^2 E[|\mathbf{s}_{k_1}^* \mathbf{s}_{k_2}|^4 \Upsilon_{k_1} \Upsilon_{k_2}] = \sum_{k_1 \neq k_2}^{K} p_{k_1}^2 p_{k_2}^2 E\{\Upsilon_{k_2} E[|\mathbf{s}_{k_1}^* \mathbf{s}_{k_2}|^4 \Upsilon_{k_1} \mid \mathbf{s}_{k_2}]\},$$

which is less than or equal to

(5.20)
$$M \sum_{k_1 \neq k_2}^{K} E\left\{|\Upsilon_{k_2}| E\left[\left(\mathbf{s}_{k_1}^* \mathbf{s}_{k_2} \mathbf{s}_{k_2}^* \mathbf{s}_{k_1} - \frac{1}{N} \operatorname{Tr} \mathbf{s}_{k_2} \mathbf{s}_{k_2}^*\right)^2 |\Upsilon_{k_1}| \mid \mathbf{s}_{k_2}\right]\right\}$$
$$+ M \sum_{k_1 \neq k_2}^{K} E\left\{|\Upsilon_{k_2}| E\left[\left(\frac{1}{N} \mathbf{s}_{k_2}^* \mathbf{s}_{k_2}\right)^2 |\Upsilon_{k_1}| \mid \mathbf{s}_{k_2}\right]\right\} = o(1),$$

as

$$E\left\{|\Upsilon_{k_2}| E\left[\left(\mathbf{s}_{k_1}^* \mathbf{s}_{k_2} \mathbf{s}_{k_2}^* \mathbf{s}_{k_1} - \frac{1}{N} \operatorname{Tr} \mathbf{s}_{k_2} \mathbf{s}_{k_2}^*\right)^2 |\Upsilon_{k_1}| \mid \mathbf{s}_{k_2}\right]\right\}$$
$$\leq E\left\{|\Upsilon_{k_2}| \left[E\left(\left(\mathbf{s}_{k_1}^* \mathbf{s}_{k_2} \mathbf{s}_{k_2}^* \mathbf{s}_{k_1} - \frac{1}{N} \operatorname{Tr} \mathbf{s}_{k_2} \mathbf{s}_{k_2}^*\right)^4 \mid \mathbf{s}_{k_2}\right)\right]^{1/2} [E(\Upsilon_{k_1}^2 \mid \mathbf{s}_{k_2})]^{1/2}\right\}$$
$$\leq \frac{M \varepsilon_N^2}{N^{3/2}} E[|\Upsilon_{k_2}| (\mathbf{s}_{k_2}^* \mathbf{s}_{k_2})^2]$$
$$\leq \frac{M \varepsilon_N^2}{N^{3/2}} E[|\Upsilon_{k_2}|^3] + \frac{M \varepsilon_N^2}{N^{3/2}} (E \Upsilon_{k_2}^2)^{1/2}$$
$$= o\left(\frac{1}{N^2}\right)$$

and

$$E\left\{|\Upsilon_{k_2}| E\left[\left(\frac{1}{N} \mathbf{s}_{k_2}^* \mathbf{s}_{k_2}\right)^2 |\Upsilon_{k_1}| \mid \mathbf{s}_{k_2}\right]\right\}$$
$$\leq \frac{1}{N^2} E[|\Upsilon_{k_2}| (\mathbf{s}_{k_2}^* \mathbf{s}_{k_2})^2] (E|\Upsilon_{k_1}|^2)^{1/2} = O\left(\frac{1}{N^3}\right).$$

Consequently, $G_{224}$ converges to zero and then $G_{22}$ converges to zero. Therefore, via (5.14),

(5.21) $$\sum_{k=1}^{K} p_k(\check{s}_k)(\mathbf{s}_k^* \mathbf{s}_k - 1) \xrightarrow{i.p.} \frac{E|v_{11}|^4 - 1}{c} a_1 \int x \, dH(x).$$

Combining (5.9)–(5.12) with (5.21), one can conclude that

$$G_2 = -\frac{1}{a_1^2} \left[ 2a_1 \frac{E|v_{11}|^4 - 1}{c} \int x \, dH(x) \right.$$

$$+ \sigma^2 \operatorname{Tr} \mathbf{C} + \operatorname{Tr} \mathbf{C}^2 - a_1 \sum_{k=1}^{K} p_k$$

(5.22)
$$- \frac{1}{c}(E|v_{11}|^4 - 1)$$

$$\times \left[ \int x^2 \, dH(x) - 2 \operatorname{Tr} \mathbf{SP}^2 \mathbf{S}^* + \sum_{k=1}^{K} p_k^2 \right] + o_p(1).$$

Fourth, turn to the term $G_3$. It is decomposed as

(5.23) $$a_1^3 G_3 = G_{31} + G_{32} + G_{33},$$

where

$$G_{31} = \sum_{k=1}^{K} p_k (\mathbf{s}_k^* \mathbf{s}_k - 1)^2 (\check{s}_k)^2$$

(recall $\check{s}_k = \mathbf{s}_k^* \mathbf{R}_k \mathbf{s}_k - a_1$) and

$$G_{32} = 2 \sum_{k=1}^{K} p_k (\mathbf{s}_k^* \mathbf{s}_k - 1)(\check{s}_k)^2, \qquad G_{33} = \sum_{k=1}^{K} p_k (\check{s}_k)^2.$$

Applying the Hölder inequality,

$$E|G_{31}| \le M \sum_{k=1}^{K} \left[ E(\mathbf{s}_k^* \mathbf{s}_k - 1)^2 \left( \mathbf{s}_k^* \mathbf{R}_k \mathbf{s}_k - \frac{1}{N} \operatorname{Tr} \mathbf{R}_k \right)^2 \right.$$

$$\left. + E(\mathbf{s}_k^* \mathbf{s}_k - 1)^2 \left( \frac{1}{N} \operatorname{Tr} \mathbf{R}_k - a_1 \right)^2 \right]$$

$$\le M \sum_{k=1}^{K} (E(\mathbf{s}_k^* \mathbf{s}_k - 1)^4)^{1/2} \left( E \left( \mathbf{s}_k^* \mathbf{R}_k \mathbf{s}_k - \frac{1}{N} \operatorname{Tr} \mathbf{R}_k \right)^4 \right)^{1/2}$$

$$+ M \sum_{k=1}^{K} E(\mathbf{s}_k^* \mathbf{s}_k - 1)^2 E \left( \frac{1}{N} \operatorname{Tr} \mathbf{R}_k - a_1 \right)^2 = o(1).$$

Analogously, one can also obtain

$$E|G_{32}| = o(1).$$

To derive the limit of $G_{33}$, we need to evaluate its variance:

(5.24) $$E \left( \sum_{k=1}^{K} p_k^2 [(\check{s}_k)^2 - E(\check{s}_k)^2] \right)^2 = G_{331} + G_{332},$$

where

$$G_{331} = \sum_{k=1}^{K} p_k^2 E[(\check{s}_k)^2 - E(\check{s}_k)^2]^2$$

and

$$G_{332} = \sum_{k_1 \neq k_2}^{K} p_{k_1} p_{k_2} E[((\check{s}_{k_1})^2 - E(\check{s}_{k_1})^2)((\check{s}_{k_2})^2 - E(\check{s}_{k_2})^2)].$$

For $G_{331}$, we have

$$G_{331} \leq M \sum_{k=1}^{K} E[(\check{s}_k)^4]$$

$$\leq M \sum_{k=1}^{K} E\left(\mathbf{s}_k^* \mathbf{R}_k \mathbf{s}_k - \frac{1}{N} \operatorname{Tr} \mathbf{R}_k\right)^4 + M \sum_{k=1}^{K} E\left(\frac{1}{N} \operatorname{Tr} \mathbf{R}_k - \frac{1}{N} \operatorname{Tr} \mathbf{R}\right)^4$$

$$+ M \sum_{k=1}^{K} E\left(\frac{1}{N} \operatorname{Tr} \mathbf{R} - a_1\right)^4 = o(1).$$

In fact, note that $a_1 = \sigma^2 + \frac{1}{c_N}$,

$$(5.25) \quad E\left(\frac{1}{N} \operatorname{Tr} \mathbf{R} - a_1\right)^4 = E\left(\frac{1}{N} \sum_{k=1}^{K} p_k (\mathbf{s}_k^* \mathbf{s}_k - 1)\right)^4 = o\left(\frac{1}{N^2}\right).$$

Since the treatment of $G_{332}$ is basically similar to that of $G_{22}$, we give only an outline. To this end, we expand it as

$$(5.26) \quad G_{332} = G_{332}^{(1)} + \cdots + G_{332}^{(9)},$$

where

$$G_{332}^{(1)} = \sum_{k_1 \neq k_2}^{K} p_{k_1} p_{k_2} \operatorname{Cov}(\alpha_{k_1}^2, \alpha_{k_2}^2),$$

$$G_{332}^{(2)} = \sum_{k_1 \neq k_2}^{K} p_{k_1}^3 p_{k_2} \operatorname{Cov}(\alpha_{k_1}^2, |\mathbf{s}_{k_1}^* \mathbf{s}_{k_2}|^4),$$

$$G_{332}^{(3)} = 2 \sum_{k_1 \neq k_2}^{K} p_{k_1}^2 p_{k_2} \operatorname{Cov}(\alpha_{k_1}^2, \alpha_{k_2} |\mathbf{s}_{k_1}^* \mathbf{s}_{k_2}|^2),$$

$$G_{332}^{(4)} = 2 \sum_{k_1 \neq k_2}^{K} p_{k_1} p_{k_2}^2 \operatorname{Cov}(\alpha_{k_1} |\mathbf{s}_{k_1}^* \mathbf{s}_{k_2}|^2, \alpha_{k_2}^2),$$

$$G_{332}^{(5)} = 2 \sum_{k_1 \neq k_2}^{K} p_{k_1}^3 p_{k_2} \operatorname{Cov}(\alpha_{k_1}|\mathbf{s}_{k_1}^*\mathbf{s}_{k_2}|^2, |\mathbf{s}_{k_1}^*\mathbf{s}_{k_2}|^4),$$

$$G_{332}^{(6)} = 4 \sum_{k_1 \neq k_2}^{K} p_{k_1}^2 p_{k_2}^2 \operatorname{Cov}(\alpha_{k_1}|\mathbf{s}_{k_1}^*\mathbf{s}_{k_2}|^2, \alpha_{k_2}|\mathbf{s}_{k_1}^*\mathbf{s}_{k_2}|^2),$$

$$G_{332}^{(7)} = \sum_{k_1 \neq k_2}^{K} p_{k_1} p_{k_2}^3 \operatorname{Cov}(|\mathbf{s}_{k_1}^*\mathbf{s}_{k_2}|^4, \alpha_{k_2}^2),$$

$$G_{332}^{(8)} = 2 \sum_{k_1 \neq k_2}^{K} p_{k_1}^3 p_{k_2}^3 \operatorname{Var}(|\mathbf{s}_{k_1}^*\mathbf{s}_{k_2}|^4)$$

and

$$G_{332}^{(9)} = 2 \sum_{k_1 \neq k_2}^{K} p_{k_1}^2 p_{k_2}^3 \operatorname{Cov}(|\mathbf{s}_{k_1}^*\mathbf{s}_{k_2}|^4, \alpha_{k_2}|\mathbf{s}_{k_1}^*\mathbf{s}_{k_2}|^2).$$

We claim that

$$G_{332} = o(1).$$

But, in the sequel, as an illustration, only terms $G_{332}^{(1)}$ and $G_{332}^{(8)}$ will be estimated, the argument for all the remaining ones are analogous and then omitted.

Using the Burkholder inequality,

$$(5.27) \quad E|\mathbf{s}_{k_1}^*\mathbf{s}_{k_2}|^8 \leq \frac{M}{N^8}\left[\left(\sum_{i=1}^{N} E|v_{ik_1}v_{ik_2}|^2\right)^4 + \sum_{i=1}^{N}(E|v_{11}|^8)^2\right]$$
$$= O\left(\frac{1}{N^3}\right),$$

which leads to

$$|G_{332}^{(8)}| \leq M \sum_{k_1 \neq k_2}^{K} p_{k_1}^2 p_{k_2}^3 E|\mathbf{s}_{k_1}^*\mathbf{s}_{k_2}|^8 = O\left(\frac{1}{N}\right).$$

From (5.15),

$$G_{332}^{(1)} = \sum_{k_1 \neq k_2}^{K} p_{k_1} p_{k_2} E[E(\alpha_{k_1}^2 - E\alpha_{k_1}^2 \mid \mathbf{R}_{k_1 k_2}) E(\alpha_{k_2}^2 - E\alpha_{k_2}^2 \mid \mathbf{R}_{k_1 k_2})]$$
$$\leq \frac{M}{N^4} \sum_{k_1 \neq k_2}^{K} \bigg[ E(\operatorname{Tr}\mathbf{R}_{k_1 k_2}^2 - E\operatorname{Tr}\mathbf{R}_{k_1 k_2}^2)^2 + E(\operatorname{Tr}\mathbf{R}_{k_1 k_2} - Na_1)^4$$

$$(5.28) \qquad + E\left(\sum_{i=1}^{N}[(\mathbf{R}_{k_1 k_2})_{ii}]^2 - E[(\mathbf{R}_{k_1 k_2})_{ii}]^2\right)^2$$

$$+ (E(\operatorname{Tr}\mathbf{R}_{k_1 k_2} - Na_1)^2)^2\Bigg]$$

$$= o(1).$$

In order to get (5.28), we need to analyze the above four terms on the right-hand side of the inequality. First, applying $\mathbf{R}_{k_1 k_2} = (\mathbf{R}_{k_1 k_2} - \mathbf{R}_{k_1}) + \mathbf{R}_{k_1}$ twice, we have

$$E(\operatorname{Tr}\mathbf{R}_{k_1 k_2}^2 - E\operatorname{Tr}\mathbf{R}_{k_1 k_2}^2)^2$$

$$\leq ME(\operatorname{Tr}(\mathbf{R}_{k_1 k_2} - \mathbf{R}_{k_1})\mathbf{R}_{k_1 k_2} - E\operatorname{Tr}(\mathbf{R}_{k_1 k_2} - \mathbf{R}_{k_1})\mathbf{R}_{k_1 k_2})^2$$

$$+ ME(\operatorname{Tr}\mathbf{R}_{k_1}(\mathbf{R}_{k_1 k_2} - \mathbf{R}_{k_1}) - E\operatorname{Tr}\mathbf{R}_{k_1}(\mathbf{R}_{k_1 k_2} - \mathbf{R}_{k_1}))^2$$

$$(5.29) \qquad + ME(\operatorname{Tr}\mathbf{R}_{k_1}^2 - E\operatorname{Tr}\mathbf{R}_{k_1}^2)^2$$

$$\leq ME(\hat{\gamma}_{k_2})^2 + ME((\mathbf{s}_{k_2}^*\mathbf{s}_{k_2})^2 - E(\mathbf{s}_{k_2}^*\mathbf{s}_{k_2})^2)^2$$

$$+ ME(\operatorname{Tr}\mathbf{R}_{k_1}^2 - E\operatorname{Tr}\mathbf{R}_{k_1}^2)^2,$$

where $\hat{\gamma}_{k_2} = \mathbf{s}_{k_2}^*\mathbf{R}_{k_1 k_2}\mathbf{s}_{k_2} - \frac{1}{N}E\operatorname{Tr}\mathbf{R}_{k_1 k_2}$. However, observe that

$$(5.30) \qquad E((\mathbf{s}_{k_2}^*\mathbf{s}_{k_2})^2 - E(\mathbf{s}_{k_2}^*\mathbf{s}_{k_2})^2)^2 = O\left(\frac{1}{N}\right)$$

and that

$$E(\gamma_{k_2})^2 \leq \frac{M}{N},$$

using (5.4) and (5.16). Therefore,

$$E(\operatorname{Tr}\mathbf{R}_{k_1 k_2}^2 - E\operatorname{Tr}\mathbf{R}_{k_1 k_2}^2)^2 \leq \frac{M}{N} + ME(\operatorname{Tr}\mathbf{R}_{k_1}^2 - E\operatorname{Tr}\mathbf{R}_{k_1}^2)^2$$

$$\leq \frac{M}{N} + ME(\operatorname{Tr}\mathbf{R}^2 - E\operatorname{Tr}\mathbf{R}^2)^2$$

again, repeating a process analogous to (5.29) in the last step. But this implies

$$\frac{M}{N^4}\sum_{k_1 \neq k_2}^{K} E(\operatorname{Tr}\mathbf{R}_{k_1 k_2}^2 - E\operatorname{Tr}\mathbf{R}_{k_1 k_2}^2)^2 \to 0.$$

Second,

$$E(\operatorname{Tr}\mathbf{R}_{k_1 k_2} - Na_1)^4$$

$$\leq ME\left|\sum_{k=1}^{K} p_k(\mathbf{s}_k^*\mathbf{s}_k - 1)\right|^4 + ME(\mathbf{s}_{k_1}^*\mathbf{s}_{k_1})^4 + ME(\mathbf{s}_{k_2}^*\mathbf{s}_{k_2})^4 \leq M,$$

which shows that the second sum in (5.28) converges to zero. Similarly, one can also prove that the fourth sum in (5.28) converges to zero by a similar argument, as expected. Finally, in order to show that the third sum in (5.28) converges to zero, it is enough to show that

$$E|[(\mathbf{C})_{ii}]^2 - E[(\mathbf{C})_{ii}]^2|^2 \to 0.$$

To this end, it suffices to verify that

$$E|(\mathbf{C})_{ii} - E(\mathbf{C})_{ii}|^4 \to 0,$$

but, as in Theorem 1.3, through martingale difference decomposition, one can get it. Thus, (5.28) holds.

Hence, by (5.15) and an argument similar to Theorem 1.4, we have so far proved that

$$G_{33} = \frac{a_2 \int x\, dH(x)}{c} + \frac{E|v_{11}|^4 - 2}{c} a_1^2 \int x\, dH(x) + o_p(1)$$

and then

(5.31) $$G_3 = \frac{1}{a_1^3 c}(a_2 + (E|v_{11}|^4 - 2)a_1^2) \int x\, dH(x) + o_p(1).$$

Summarizing (5.7), (5.22) and (5.31), we conclude that

(5.32)
$$\sum_{k=1}^{K}\left(\beta_k - \frac{p_k}{\alpha_1}\right)$$
$$= \left(\frac{2}{a_1} - \frac{\sigma^2}{a_1^2}\right)\mathrm{Tr}\,\mathbf{C} - \frac{1}{a_1^2}\mathrm{Tr}\,\mathbf{C}^2 + \frac{2}{a_1^2}\mathrm{Tr}\,\mathbf{SP}^2\mathbf{S}^*$$
$$- \frac{1}{a_1}\sum_{k=1}^{K} p_k - \frac{1}{a_1^2}\sum_{k=1}^{K} p_k^2 + \left(\frac{a_2}{ca_1^3} - \frac{1}{ca_1}\right)\int x\, dH(x)$$
$$+ \frac{E|v_{11}|^4 - 1}{ca_1^2}\int x^2\, dH(x) + o_p(1)$$

(recall that $a_1 = \sigma^2 + 1/c$).

Now we let $p_k = 1, k = 1, \ldots, K$, so that Theorem 1.4 can be applied. In this case (5.32) becomes

(5.33)
$$\sum_{k=1}^{K}\left(\beta_k - \frac{p_k}{\alpha_1}\right) = \frac{2 + 2/c + \sigma^2}{a_1^2}(\mathrm{Tr}\,\mathbf{C} - K) - \frac{1}{a_1^2}(\mathrm{Tr}\,\mathbf{C}^2 - (1 + 1/c)K)$$
$$+ \frac{1}{c^2 a_1^3} + \frac{E|v_{11}|^4 - 1}{ca_1^2} + o_p(1)$$

$[a_2 = \frac{(1+1/c)}{c} + \sigma^4 + 2\sigma^2/c]$. Then Theorem 1.3 follows from Theorem 1.4. The expectation and variance of the limiting normal distribution are also from (1.23) and (1.24) of Bai and Silverstein [3] and our (1.23) and (1.24). Therefore, the variance $\tau^2$ is equal to

$$
\begin{aligned}
&\frac{(2+2/c+\sigma^2)^2(E|v_{11}|^4 - 1)}{c(\sigma^2+1/c)^4} \\
&\quad - 2\frac{(2+2/c+\sigma^2)(2c+2)}{c^2(\sigma^2+1/c)^4}(E|v_{11}|^4 - 1) \\
&\quad + \frac{1}{c^4(\sigma^2+1/c)^4}(4c^3 + 10c^2 + 4c \\
&\qquad\qquad + (4c^3 + 8c^2 + 4c)(E|v_{11}|^4 - 2)).
\end{aligned}
\tag{5.34}
$$

**6. Proof of Corollary 1.1.** By the Taylor expansion,

$$
\sum_{k=1}^K \left(\log(1+\beta_k) - \log\left(1 + \frac{1}{a_1}\right)\right)
$$

$$
= \sum_{k=1}^K \frac{(\beta_k - 1/a_1)}{1 + 1/a_1} - \sum_{k=1}^K \frac{(\beta_k - 1/a_1)^2}{2(1+1/a_1)^2} + \sum_{k=1}^K \frac{(\beta_k - 1/a_1)^3}{3(1+\psi_k)^3},
$$

with $\psi_k$ being located in the interval $[1/a_1, \beta_k]$. Similar to Theorem 1.3, it can be shown that

$$
E\left(\sum_{k=1}^K \left(\left(\beta_k - \frac{1}{a_1}\right)^2 - E\left(\beta_k - \frac{1}{a_1}\right)^2\right)\right)^2 = o(1)
$$

and

$$
\sum_{k=1}^K \frac{|(\beta_k - 1/a_1)^3|}{3(1+\psi_k)^3} \leq M \sum_{k=1}^K \left|\left(\beta_k - \frac{1}{a_1}\right)^3\right| \xrightarrow{i.p.} 0.
$$

Now compute $\sum_{k=1}^K E(\beta_k - \frac{1}{a_1})^2$. Applying (5.1),

$$
(6.1) \quad \sum_{k=1}^K E\left(\beta_k - \frac{1}{a_1}\right)^2 = \sum_{k=1}^K E\left(\frac{(\mathbf{s}_k^*\mathbf{s}_k)^2 - 1}{a_1} - \frac{(\mathbf{s}_k^*\mathbf{s}_k)^2(\check{s}_k)}{a_1^2} + \frac{(\mathbf{s}_k^*\mathbf{s}_k)^2(\check{s}_k)^2}{a_1^2 \mathbf{s}_k^* \mathbf{R}_k \mathbf{s}_k}\right)^2.
$$

Again, via an argument analogous to Theorem 1.3, it is easily seen that the contribution from the terms involving $(\mathbf{s}_k^*\mathbf{s}_k)^2(\check{s}_k)^2/(a_1^2 \mathbf{s}_k^*\mathbf{R}_k\mathbf{s}_k)$ can be ignored. So (6.1) is equal to

$$
\sum_{k=1}^K E\left(\frac{(\mathbf{s}_k^*\mathbf{s}_k)^2 - 1}{a_1}\right)^2 + \sum_{k=1}^K E\left(\frac{(\mathbf{s}_k^*\mathbf{s}_k)^2(\check{s}_k)}{a_1^2}\right)^2
$$

$$- 2\sum_{k=1}^{K} E \frac{(\mathbf{s}_k^* \mathbf{s}_k)^2 - 1}{a_1} \frac{(\mathbf{s}_k^* \mathbf{s}_k)^2 (\check{s}_k)}{a_1^2} + o(1).$$

Combining the argument of Theorem 1.3 with (5.15) and (5.14), the above three terms are equal to, respectively, $4(E|v_{11}|^4 - 1)/(ca_1^2)$, $(a_2 + E|v_{11}|^4 - 2)/(ca_1^4)$ and $-4(E|v_{11}|^4 - 1)/(ca_1^2)$. Thus, we finish the proof.

## APPENDIX: TRUNCATION OF UNDERLYING RANDOM VARIABLES IN THEOREM 1.3

Let $\hat{v}_{ij} = v_{ij} I(|v_{ij}| \leq \varepsilon_N \sqrt{N})$ and $\bar{v}_{ij} = \hat{v}_{ij} - E\hat{v}_{ij}$, $i = 1, \ldots, N$, $j = 1, \ldots, K$, where $\varepsilon_N$ is a positive sequence converging to zero. We use $\hat{\mathbf{s}}_k, \bar{\mathbf{s}}_k, \hat{\mathbf{S}}_k, \bar{\mathbf{S}}_k, \hat{\mathbf{R}}_k, \bar{\mathbf{R}}_k$ and $\hat{\beta}_k, \bar{\beta}_k$, $k = 1, \ldots, K$, to denote the analogues of $\mathbf{s}_k$, $\mathbf{S}_k$, $\mathbf{R}_k$ and $\beta_k$ with the elements replaced by $\hat{v}_{ij}$ or $\bar{v}_{ij}$.

As in the proof for Theorem 1.3 in Pan, Guo and Zhou [12], one can select the above $\varepsilon_N$ so that

$$\varepsilon_N^{-4} E v_{11}^4 I(|v_{11}| > \varepsilon_N \sqrt{N}) \to 0,$$

and show that

$$\sum_{k=1}^{K} \beta_k - \sum_{k=1}^{K} \hat{\beta}_k \xrightarrow{i.p.} 0.$$

Now consider the re-centralization of random variables. Applying (5.1),

(A.1) $$\sum_{k=1}^{K} \frac{p_k (\hat{\mathbf{s}}_k^* \hat{\mathbf{s}}_k)^2}{\hat{\mathbf{s}}_k^* \hat{\mathbf{R}}_k \hat{\mathbf{s}}_k} = U_1 + U_2 + U_3 + U_4,$$

where

$$U_1 = \frac{1}{a_1} \sum_{k=1}^{K} p_k (\hat{\mathbf{s}}_k^* \hat{\mathbf{s}}_k)^2, \qquad U_4 = -\frac{1}{a_1^3} \sum_{k=1}^{K} \frac{p_k (\hat{\mathbf{s}}_k^* \hat{\mathbf{s}}_k)^2 (\hat{\mathbf{s}}_k^* \hat{\mathbf{R}}_k \hat{\mathbf{s}}_k - a_1)^3}{\hat{\mathbf{s}}_k^* \hat{\mathbf{R}}_k \hat{\mathbf{s}}_k}$$

and

$$U_2 = -\frac{1}{a_1^2} \sum_{k=1}^{K} p_k (\hat{\mathbf{s}}_k^* \hat{\mathbf{s}}_k)^2 (\hat{\mathbf{s}}_k^* \hat{\mathbf{R}}_k \hat{\mathbf{s}}_k - a_1),$$

$$U_3 = \frac{1}{a_1^3} \sum_{k=1}^{K} p_k (\hat{\mathbf{s}}_k^* \hat{\mathbf{s}}_k)^2 (\hat{\mathbf{s}}_k^* \hat{\mathbf{R}}_k \hat{\mathbf{s}}_k - a_1)^2.$$

In the sequel we shall show that $U_4$ converges to zero in probability. Note that

$$\frac{\hat{\mathbf{s}}_k^* \hat{\mathbf{s}}_k}{\hat{\mathbf{s}}_k^* \hat{\mathbf{R}}_k \hat{\mathbf{s}}_k} \leq \frac{1}{\sigma^2}, \qquad \hat{\mathbf{s}}_k^* = \bar{\mathbf{s}}_k^* + E\hat{\mathbf{s}}_k^*.$$

This gives

$$|U_4| \leq M \sum_{k=1}^{K} \hat{\mathbf{s}}_k^* \hat{\mathbf{s}}_k |\hat{\mathbf{s}}_k^* \hat{\mathbf{R}}_k \hat{\mathbf{s}}_k - a_1|^3$$

$$(A.2) \qquad \leq M \sum_{k=1}^{K} (\bar{\mathbf{s}}_k^* \bar{\mathbf{s}}_k + \bar{\mathbf{s}}_k^* E\hat{\mathbf{s}}_k + (E\hat{\mathbf{s}}_k^*)\bar{\mathbf{s}}_k + E\hat{\mathbf{s}}_k^* E\hat{\mathbf{s}}_k)$$

$$\times |\bar{\mathbf{s}}_k^* \hat{\mathbf{R}}_k \bar{\mathbf{s}}_k - a_1 + \bar{\mathbf{s}}_k^* \hat{\mathbf{R}}_k E\hat{\mathbf{s}}_k + (E\hat{\mathbf{s}}_k^*)\hat{\mathbf{R}}_k \bar{\mathbf{s}}_k + (E\hat{\mathbf{s}}_k^*)\hat{\mathbf{R}}_k E\hat{\mathbf{s}}_k|^3.$$

Then we need to compute each term of the above expansion.

It is observed that

$$(A.3) \quad \left( \sum_{k=1}^{K} \bar{\mathbf{s}}_k^* \bar{\mathbf{s}}_k \left| \bar{\mathbf{s}}_k^* \hat{\mathbf{R}}_k \bar{\mathbf{s}}_k - \frac{1}{N} \operatorname{Tr} \hat{\mathbf{R}}_k \right|^3 \right)$$

$$\leq M \left( \sum_{k=1}^{K} \bar{\mathbf{s}}_k^* \bar{\mathbf{s}}_k \left( \left| \bar{\mathbf{s}}_k^* \bar{\mathbf{R}}_k \bar{\mathbf{s}}_k - \frac{1}{N} \operatorname{Tr} \bar{\mathbf{R}}_k \right|^3 + |\bar{\mathbf{s}}_k^* (E\hat{\mathbf{S}}_k) \mathbf{P}_k \bar{\mathbf{S}}_k^* \bar{\mathbf{s}}_k|^3 \right. \right.$$

$$(A.4) \qquad \left. \left. + |\bar{\mathbf{s}}_k^* \bar{\mathbf{S}}_k \mathbf{P}_k (E\hat{\mathbf{S}}_k^*) \bar{\mathbf{s}}_k|^3 + |\bar{\mathbf{s}}_k^* (E\hat{\mathbf{S}}_k) \mathbf{P}_k (E\hat{\mathbf{S}}_k^*) \bar{\mathbf{s}}_k|^3 \right) \right).$$

It is a simple matter to prove that

$$\lim_{N \to \infty} E(\bar{v}_{11})^2 = 1.$$

Then, appealing (5.4), we have

$$(A.5) \quad E \left( \sum_{k=1}^{K} \bar{\mathbf{s}}_k^* \bar{\mathbf{s}}_k \left| \bar{\mathbf{s}}_k^* \bar{\mathbf{R}}_k \bar{\mathbf{s}}_k - \frac{1}{N} \operatorname{Tr} \bar{\mathbf{R}}_k \right|^3 \right)$$

$$\leq M \sum_{k=1}^{K} (E|\bar{\mathbf{s}}_k^* \bar{\mathbf{s}}_k - 1|^2)^{1/2} \left( E \left| \bar{\mathbf{s}}_k^* \bar{\mathbf{R}}_k \bar{\mathbf{s}}_k - \frac{1}{N} \operatorname{Tr} \bar{\mathbf{R}}_k \right|^6 \right)^{1/2} = o\left( \frac{1}{\sqrt{N}} \right).$$

For the first term in (A.4), with the notation $\mathbf{e} = \frac{1}{\sqrt{N}}(1, \ldots, 1)^*$ and $\mathbf{G} = \bar{\mathbf{S}}_k \mathbf{P}_k (E\hat{\mathbf{S}}_k^*)(E\hat{\mathbf{S}}_k) \mathbf{P}_k \bar{\mathbf{S}}_k^*$, analogously,

$$E \left( \sum_{k=1}^{K} \bar{\mathbf{s}}_k^* \bar{\mathbf{s}}_k |\bar{\mathbf{s}}_k^* \bar{\mathbf{S}}_k \mathbf{P}_k (E\hat{\mathbf{S}}_k^*) \bar{\mathbf{s}}_k|^3 \right)$$

$$\leq M E \left( \sum_{k=1}^{K} \bar{\mathbf{s}}_k^* \bar{\mathbf{s}}_k \left| \bar{\mathbf{s}}_k^* \bar{\mathbf{S}}_k \mathbf{P}_k (E\hat{\mathbf{S}}_k^*) \bar{\mathbf{s}}_k - \frac{1}{N} \operatorname{Tr} \bar{\mathbf{S}}_k \mathbf{P}_k E\hat{\mathbf{S}}_k^* \right|^3 \right)$$

$$+ M E \left( \sum_{k=1}^{K} \bar{\mathbf{s}}_k^* \bar{\mathbf{s}}_k \left| \frac{1}{N} \operatorname{Tr} \bar{\mathbf{S}}_k \mathbf{P}_k E\hat{\mathbf{S}}_k^* \right|^3 \right)$$

$$\leq \frac{M}{N^{3/2}} \sum_{k=1}^{K} \left[ \left( E\left(\frac{1}{N}\operatorname{Tr}\mathbf{G}\right)^3 + E\frac{1}{N}\operatorname{Tr}(\mathbf{G})^3 \right)^{1/2} \right.$$

$$\left. + E\left(\frac{1}{N}\operatorname{Tr}\mathbf{G}\right)^{3/2} + E\frac{1}{N}\operatorname{Tr}(\mathbf{G})^{3/2} \right]$$

$$+ M\frac{(E\hat{v}_{11})^3}{N^3}\left(\sum_{k=1}^{K} E\left|\sum_{j\neq k} p_j \bar{s}_j^* e\right|^3\right) = o(1).$$

Indeed,

$$(E\hat{v}_{11})^3 \left(\sum_{k=1}^{K} E\left|\sum_{j\neq k} p_j \bar{s}_j^* e\right|^3\right) \leq M(E\hat{v}_{11})^3 K N^{3/2} = o\left(\frac{1}{N^2}\right),$$

and since

$$\mathbf{G} \leq MN(E\hat{v}_{11})^2 \bar{\mathbf{S}}_k \mathbf{e}\mathbf{e}^* \mathbf{e}\mathbf{e}^* \bar{\mathbf{S}}_k^*,$$

we have

$$E\left(\frac{1}{N}\operatorname{Tr}\mathbf{G}\right)^3 \leq MN^3(E\hat{v}_{11})^6 E\left(\frac{1}{N}\operatorname{Tr}\bar{\mathbf{S}}_1\bar{\mathbf{S}}_1^*\right)^3.$$

The other terms can be estimated similarly. Regarding the last term in (A.4), we have

$$\sum_{k=1}^{K} \bar{\mathbf{s}}_k^* \bar{\mathbf{s}}_k |\bar{\mathbf{s}}_k^*(E\hat{\mathbf{S}}_k)\mathbf{P}_k(E\hat{\mathbf{S}}_k^*)\bar{\mathbf{s}}_k|^3 \leq M \sum_{k=1}^{K} (\bar{\mathbf{s}}_k^*\bar{\mathbf{s}}_k)^4 \|E\hat{\mathbf{S}}_k\|^6 = o\left(\frac{1}{n}\right),$$

where $\|\cdot\|$ denotes the spectral norm of a matrix. Therefore, (A.3) converges to zero in probability.

Now

$$\sum_{k=1}^{K} \bar{\mathbf{s}}_k^* \bar{\mathbf{s}}_k \left|\frac{1}{N}\operatorname{Tr}\hat{R}_k - a_1\right|^3$$

(A.6)
$$\leq M \sum_{k=1}^{K} \bar{\mathbf{s}}_k^* \bar{\mathbf{s}}_k \left( \left|\frac{1}{N}\operatorname{Tr}\bar{R}_k - a_1\right|^3 + \left|\frac{1}{N}\operatorname{Tr}\bar{\mathbf{S}}_k \mathbf{P}_k E\hat{\mathbf{S}}_k^*\right|^3 \right.$$

$$\left. + \left|\frac{1}{N}\operatorname{Tr}(E\hat{\mathbf{S}}_k)\mathbf{P}_k\bar{\mathbf{S}}_k^*\right|^3 + \left|\frac{1}{N}\operatorname{Tr}E\hat{\mathbf{S}}_k\mathbf{P}_k E\hat{\mathbf{S}}_k^*\right|^3 \right).$$

For its first term one can get

$$E\sum_{k=1}^{K} \bar{\mathbf{s}}_k^* \bar{\mathbf{s}}_k \left|\frac{1}{N}\operatorname{Tr}\bar{R}_k - a_1\right|^3 \leq M \sum_{k=1}^{K} (E(\bar{\mathbf{s}}_k^*\bar{\mathbf{s}}_k - 1)^2)^{1/2} \left(E\left|\frac{1}{N}\operatorname{Tr}\bar{R}_k - a_1\right|^6\right)^{1/2}$$

$$+ M\sum_{k=1}^{K} E\left|\frac{1}{N}\operatorname{Tr}\bar{R}_k - a_1\right|^3 = o(1),$$

as

$$E\left|\frac{1}{N}\operatorname{Tr}\bar{R}_k - a_1\right|^p = E\left|\frac{1}{N}\sum_{j\neq k}^K p_j\bar{\mathbf{s}}_j^*\bar{\mathbf{s}}_j - c\right|^p = O\left(\frac{1}{N^{p/2}}\right).$$

The argument that the remaining terms of (A.6) converge to zero in probability is similar to above, even simpler and then omitted. Hence, (A.6) converges to zero in probability. This, together with (A.3), leads to

$$\text{(A.7)} \qquad \sum_{k=1}^K \bar{\mathbf{s}}_k^*\bar{\mathbf{s}}_k |\bar{\mathbf{s}}_k^*\hat{\mathbf{R}}_k\bar{\mathbf{s}}_k - a_1|^3 \xrightarrow{i.p.} 0,$$

which is one term in (A.2). All remaining items of (A.2) can be computed similarly, so we omit it here. Consequently,

$$U_4 \xrightarrow{i.p.} 0$$

and

$$\text{(A.8)} \qquad \sum_{k=1}^K \frac{p_k(\hat{\mathbf{s}}_k^*\hat{\mathbf{s}}_k)^2}{\hat{\mathbf{s}}_k^*\hat{\mathbf{R}}_k\hat{\mathbf{s}}_k} = U_1 + U_2 + U_3 + o_p(1).$$

Analogously, one can also show that

$$\text{(A.9)} \qquad \sum_{k=1}^K \frac{p_k(\bar{\mathbf{s}}_k^*\bar{\mathbf{s}}_k)^2}{\bar{\mathbf{s}}_k^*\bar{\mathbf{R}}_k\bar{\mathbf{s}}_k} = V_1 + V_2 + V_3 + o_p(1),$$

where

$$V_1 = \frac{1}{a_1}\sum_{k=1}^K p_k(\bar{\mathbf{s}}_k^*\bar{\mathbf{s}}_k)^2, \qquad V_2 = -\frac{1}{a_1^2}\sum_{k=1}^K p_k(\bar{\mathbf{s}}_k^*\bar{\mathbf{s}}_k)^2(\bar{\mathbf{s}}_k^*\bar{\mathbf{R}}_k\bar{\mathbf{s}}_k - a_1)$$

and

$$V_3 = \frac{1}{a_1^3}\sum_{k=1}^K p_k(\bar{\mathbf{s}}_k^*\bar{\mathbf{s}}_k)^2(\bar{\mathbf{s}}_k^*\bar{\mathbf{R}}_k\bar{\mathbf{s}}_k - a_1)^2.$$

In the following, we show that $U_i - V_i, i = 1, 2, 3$, converge to zero in probability. Since all the calculations for $U_i - V_i$ are similar, as an illustration, we consider $U_2 - V_2$ only.

Write

$$\text{(A.10)} \qquad -a_1^2(U_2 - V_2) = U_{21} + U_{22},$$

where

$$U_{21} = \sum_{k=1}^K p_k((\hat{\mathbf{s}}_k^*\hat{\mathbf{s}}_k)^2 - (\bar{\mathbf{s}}_k^*\bar{\mathbf{s}}_k)^2)(\bar{\mathbf{s}}_k^*\bar{\mathbf{R}}_k\bar{\mathbf{s}}_k - a_1)$$

and

$$U_{22} = \sum_{k=1}^{K} p_k (\hat{\mathbf{s}}_k^* \hat{\mathbf{s}}_k)^2 (\hat{\mathbf{s}}_k^* \hat{\mathbf{R}}_k \hat{\mathbf{s}}_k - \bar{\mathbf{s}}_k^* \bar{\mathbf{R}}_k \bar{\mathbf{s}}_k).$$

Furthermore, expand $U_{21}$ as follows:

$$\begin{aligned} U_{21} = \sum_{k=1}^{K} p_k (&(\bar{\mathbf{s}}_k^* E \hat{\mathbf{s}}_k)^2 + ((E \hat{\mathbf{s}}_k^*) \bar{\mathbf{s}}_k)^2 + (E \hat{\mathbf{s}}_k^* E \hat{\mathbf{s}}_k)^2 + 2 \bar{\mathbf{s}}_k^* \bar{\mathbf{s}}_k \bar{\mathbf{s}}_k^* E \hat{\mathbf{s}}_k \\ &+ 2 \bar{\mathbf{s}}_k^* \bar{\mathbf{s}}_k (E \hat{\mathbf{s}}_k^*) \bar{\mathbf{s}}_k + 2 \bar{\mathbf{s}}_k^* (E \hat{\mathbf{s}}_k)(E \hat{\mathbf{s}}_k^*) \bar{\mathbf{s}}_k + 2 \bar{\mathbf{s}}_k^* \bar{\mathbf{s}}_k E \hat{\mathbf{s}}_k^* E \hat{\mathbf{s}}_k \\ &+ 2 (E \hat{\mathbf{s}}_k^*) \bar{\mathbf{s}}_k (E \hat{\mathbf{s}}_k^*) E \hat{\mathbf{s}}_k + 2 \bar{\mathbf{s}}_k^* E \hat{\mathbf{s}}_k E \hat{\mathbf{s}}_k^* E \hat{\mathbf{s}}_k) \\ \times\ & (\bar{\mathbf{s}}_k^* \bar{\mathbf{R}}_k \bar{\mathbf{s}}_k - a_1), \end{aligned}$$

which can be easily proved to tend to zero in probability. For example, for one of the terms,

$$\begin{aligned} E \left| \sum_{k=1}^{K} p_k \bar{\mathbf{s}}_k^* \bar{\mathbf{s}}_k \bar{\mathbf{s}}_k^* (E \hat{\mathbf{s}}_k)(\bar{\mathbf{s}}_k^* \bar{\mathbf{R}}_k \bar{\mathbf{s}}_k - a_1) \right| \\ \leq \sqrt{N} \varepsilon_N E \hat{v}_{11} \sum_{k=1}^{K} (E(\bar{\mathbf{s}}_k^* \bar{\mathbf{s}}_k)^2)^{1/2} (E(\bar{\mathbf{s}}_k^* \bar{\mathbf{R}}_k \bar{\mathbf{s}}_k - a_1)^2)^{1/2} = o\left(\frac{1}{\sqrt{N}}\right). \end{aligned}$$

For $U_{22}$, we have

(A.11)
$$\begin{aligned} U_{22} \leq \sum_{k=1}^{K} p_k (\hat{\mathbf{s}}_k^* \hat{\mathbf{s}}_k)^2 (&\bar{\mathbf{s}}_k^* \hat{\mathbf{R}}_k \bar{\mathbf{s}}_k - \bar{\mathbf{s}}_k^* \bar{\mathbf{R}}_k \bar{\mathbf{s}}_k + \bar{\mathbf{s}}_k^* \hat{\mathbf{R}}_k E \hat{\mathbf{s}}_k \\ &+ (E \hat{\mathbf{s}}_k^*) \hat{\mathbf{R}}_k \bar{\mathbf{s}}_k + (E \bar{\mathbf{s}}_k^*) \hat{\mathbf{R}}_k (E \hat{\mathbf{s}}_k)). \end{aligned}$$

Moreover, it is observed that

$$\sum_{k=1}^{K} p_k (\hat{\mathbf{s}}_k^* \hat{\mathbf{s}}_k)^2 (\bar{\mathbf{s}}_k^* \hat{\mathbf{R}}_k \bar{\mathbf{s}}_k - \bar{\mathbf{s}}_k^* \bar{\mathbf{R}}_k \bar{\mathbf{s}}_k)$$

(A.12)
$$= \sum_{k=1}^{K} p_k (\hat{\mathbf{s}}_k^* \hat{\mathbf{s}}_k)^2$$
$$\times (\bar{\mathbf{s}}_k^* \bar{\mathbf{S}}_k \mathbf{P}_k (E \hat{\mathbf{S}}_k^*) \bar{\mathbf{s}}_k + \bar{\mathbf{s}}_k^* (E \hat{\mathbf{S}}_k) \mathbf{P}_k \bar{\mathbf{S}}_k^* \bar{\mathbf{s}}_k + \bar{\mathbf{s}}_k^* (E \hat{\mathbf{S}}_k) \mathbf{P}_k (E \hat{\mathbf{S}}_k^*) \bar{\mathbf{s}}_k).$$

As for the first sum of the above expansion,

$$E \left| \sum_{k=1}^{K} p_k (\hat{\mathbf{s}}_k^* \hat{\mathbf{s}}_k)^2 \bar{\mathbf{s}}_k^* \bar{\mathbf{S}}_k \mathbf{P}_k (E \hat{\mathbf{S}}_k^*) \bar{\mathbf{s}}_k \right|$$

$$\leq M \sum_{k=1}^{K} (E(\hat{\mathbf{s}}_k^* \hat{\mathbf{s}}_k)^4)^{1/2} \bigg( E \bigg( \bar{\mathbf{s}}_k^* \bar{\mathbf{S}}_k \mathbf{P}_k (E\hat{\mathbf{S}}_k^*) \bar{\mathbf{s}}_k - \frac{1}{N} \operatorname{Tr} \bar{\mathbf{S}}_k \mathbf{P}_k E\hat{\mathbf{S}}_k^* \bigg)^2 \bigg)^{1/2}$$

$$+ E \bigg| \sum_{k=1}^{K} p_k (\hat{\mathbf{s}}_k^* \hat{\mathbf{s}}_k)^2 \frac{1}{N} \operatorname{Tr} \bar{\mathbf{S}}_k \mathbf{P}_k E\hat{\mathbf{S}}_k^* \bigg| = O\bigg(\frac{1}{N^2}\bigg),$$

where we make use of

$$E(\hat{\mathbf{s}}_k^* \hat{\mathbf{s}}_k)^4 \leq M E(\hat{\mathbf{s}}_k^* \hat{\mathbf{s}}_k - E(\hat{\mathbf{s}}_k^* \hat{\mathbf{s}}_k))^4 + M(E(\hat{\mathbf{s}}_k^* \hat{\mathbf{s}}_k))^4 \leq \frac{M}{N} + M$$

and

$$E\bigg( \bar{\mathbf{s}}_k^* \bar{\mathbf{S}}_k \mathbf{P}_k (E\hat{\mathbf{S}}_k^*) \bar{\mathbf{s}}_k - \frac{1}{N} \operatorname{Tr} \bar{\mathbf{S}}_k \mathbf{P}_k E\hat{\mathbf{S}}_k^* \bigg)^2 \leq \frac{M}{N^2} \operatorname{Tr} E\hat{\mathbf{S}}_k^* E\hat{\mathbf{S}}_k E(\bar{\mathbf{S}}_k^* \bar{\mathbf{S}}_k)$$

$$\leq \frac{M \|E\hat{\mathbf{S}}_k\|^2}{N^2} E \operatorname{Tr} \bar{\mathbf{S}}_k \bar{\mathbf{S}}_k = o\bigg(\frac{1}{N^3}\bigg)$$

and

$$E \bigg| \sum_{k=1}^{K} p_k (\hat{\mathbf{s}}_k^* \hat{\mathbf{s}}_k)^2 \frac{1}{N} \operatorname{Tr} \bar{\mathbf{S}}_k \mathbf{P}_k E\hat{\mathbf{S}}_k^* \bigg| = o\bigg(\frac{1}{N^{1/2}}\bigg).$$

Similarly, one can also verify that the other two terms of (A.12) converge to zero in probability, and all the other items of (A.11) converge to zero in probability. So $U_{22} \xrightarrow{i.p.} 0$, and then $U_2 - V_2 \xrightarrow{i.p.} 0$, as expected. Finally, we get

$$\sum_{k=1}^{K} \hat{\beta}_k - \sum_{k=1}^{K} \bar{\beta}_k \xrightarrow{i.p.} 0.$$

Similarly, one can perform the re-normalization step, but it is omitted here.

**Acknowledgments.** The authors would like to thank a referee for his valuable comments, which greatly improved the exposition of this work.

EURANDOM  
P.O. BOX 513  
5600MB EINDHOVEN  
THE NETHERLANDS  
E-MAIL: stapgm@gmail.com

DEPARTMENT OF STATISTICS  
AND APPLIED PROBABILITY  
NATIONAL UNIVERSITY OF SINGAPORE  
SINGAPORE 117546  
E-MAIL: stazw@nus.edu.sg